\documentclass[letterpaper,10pt]{article}

\usepackage[english]{babel}
\usepackage{amssymb}
\usepackage{latexsym}
\usepackage{amsfonts}
\usepackage{amsmath}
\usepackage{amscd}
\usepackage[dvips]{graphicx}

\newcommand{\FAT}[1]{\mathbb{#1}}
\newcommand{\hsm}[1]{\mathcal{#1}}
\newcommand{\CC}{\FAT{C}}
\newcommand{\EE}{\FAT{E}}
\newcommand{\NN}{\FAT{N}}
\newcommand{\ZZ}{\FAT{Z}}
\newcommand{\RR}{\FAT{R}}
\newcommand{\HH}{\hsm{H}}
\renewcommand{\SS}{\FAT{S}}
\newcommand{\THEN}{\,\Rightarrow\,}
\newcommand{\IFF}{\,\Leftrightarrow\,}

\newcommand{\ep}{\mbox{$\quad\blacksquare$}}

\newcommand{\proof}[1]{\textbf{Proof #1: }}
\newcommand{\minus}{\smallsetminus}
\newcommand{\cl}[1]{\overline{#1}}
\newcommand{\bd}{\partial_\infty}

\renewcommand{\ae}{\stackrel{ae}{=}}

\newcommand{\med}{\mathbf{med}}
\newcommand{\sep}[3]{\mbox{$\,#1\ast #2\ast #3\,$}}

\newcommand{\sh}[1]{\mathrm{sh}\left(#1\right)}

\newcommand{\shcirc}[1]{\mathrm{sh}^\circ\left(#1\right)}

\newtheorem{thm}{Theorem}[section]

\newtheorem{question}{Question}[section]
\newtheorem{prop}[thm]{Proposition}
\newtheorem{lemma}[thm]{Lemma}
\newtheorem{defn}[thm]{Definition}
\newtheorem{example}[thm]{Example}
\newtheorem{cor}[thm]{Corollary}
\newtheorem{remark}[thm]{Remark}

\begin{document}
\title{Coarse decompositions for boundaries of CAT(0) groups}
\author{Dan P. Guralnik\\ Vanderbilt University}
\maketitle
\begin{abstract} In this work we introduce a new combinatorial notion of boundary $\Re C$ of an $\omega$-dimensional cubing $C$. $\Re C$ is defined to be the set of almost-equality classes of ultrafilters on the standard system of halfspaces of $C$, endowed with an order relation reflecting the interaction between the Tychonoff closures of the classes.

When $C$ arises as the dual of a cubulation -- or discrete system of halfspaces -- $\HH$ of a CAT(0) space $X$ (for example, the Niblo-Reeves cubulation of the Davis-Moussong complex of a finite rank Coxeter group), we show how $\HH$ induces a function $\rho:\bd X\to\Re C$. We develop a notion of uniformness for $\HH$, generalizing the parallel walls property enjoyed by Coxeter groups, and show that, if the pair $(X,\HH)$ admits a geometric action by a group $G$, then the fibers of $\rho$ form a stratification of $\bd X$ graded by the order structure of $\Re C$. We also show how this structure computes the components of the Tits boundary of $X$.

Finally, using our result from another paper, that the uniformness of a cubulation as above implies the local finiteness of $C$, we give a condition for the co-compactness of the action of $G$ on $C$ in terms of $\rho$, generalizing a result of Williams, previously known only for Coxeter groups.
\end{abstract}

\section{Introduction}
The current research project began with the aim of exploring the relation between the asymptotic topology of a non-positively curved group $G$ and its ability to split (as an amalgam or an HNN extension) over a finitely generated quasi-convex subgroup. A group $G$ is said to be non-positively curved, if it acts properly and co-compactly by isometries (i.e., -- \emph{geometrically}) on a CAT(0) space $X$. The class of CAT(0) groups may be regarded (though the level to which this is true remains an important open question in geometric group theory) as a generalization of the class of (strongly) relatively-hyperbolic groups.

For this reason, of particular interest are splittings of CAT(0) groups over finitely generated virtually-abelian subgroups: a hyperbolic group cannot contain a free abelian subgroup of rank $2$, while a relatively-hyperbolic group may only contain such a subgroup within a parabolic subgroup. This fact has a serious impact on the connectivity properties of the boundaries of such groups (see \cite{[Bo1],[Gu]}). In this respect CAT(0) groups are more flexible, in the sense that they admit abelian subgroups of rank greater than one, though in a controlled fashion: any subgroup $H\cong \ZZ^d$ in a group $G$ acting geometrically on a CAT(0) space $X$ necessarily stabilizes an isometrically-embedded copy $F$ of $\EE^d$, on which it acts by translations so that $F/H$ is isometric to a $d$-torus.

Another reason for considering finitely generated abelian subgroups as candidates for splittings lies in the important works of Rips and Sela, where it is shown that every finitely presented torsion-free group (initially, every torsion-free word-hyperbolic group) has a canonical decomposition as a graph of groups with virtually abelian edge groups (the ``JSJ splitting'' of the group), which is, in some sense, maximal among all such decompositions. This fact lies at the base of Sela's approach to constructing algebraic geometry over groups.\\

Let $G$ denote a finitely generated group and let $\Gamma$ be a Cayley graph of $G$ with respect to a fixed finite generating set. Note that the natural action of $G$ on $\Gamma$ by left translations is a geometric action. Suppose now $X$ is a geodesic metric space admitting a geometric action by $G$; it is then basic to the approach of geometric group theory to identify $\Gamma$ with $X$, because, by the \v{S}varc-Milnor lemma, $\Gamma$ is quasi-isometric to $X$. As a result, many properties of $X$, including topological properties of compactifications of $X$, influence the asymptotics of $\Gamma$. For example, $G$ is word-hyperbolic iff $X$ is a Gromov-hyperbolic metric space, and the Gromov boundary of $X$ turns out to be equivariantly homeomorphic to the Gromov boundary of $G$.\\

The asymptotic topological theory of splittings over finite subgroups turned out to be rather simple, eventually. Given our group $G$ as above, one considers the space of ends of $\Gamma$, or, rather, its cardinality $e(G)$: a classical argument by Hopf quickly yields the following classification:
\begin{description}
    \item[$e(G)=0$.] $G$ is a finite group;
    \item[$e(G)=1$.] $G$ is infinite, any compactification of $G$ has connected boundary;
    \item[$e(G)=2$.] $G$ is virtually-cyclic;
    \item[$e(G)=\infty$.] Actually, $e(G)=2^{\aleph_0}$, as the endspace of $\Gamma$ can be shown to be perfect.
\end{description}

In the latter two cases, Stallings' theorem tells us that $G$ splits over a finite subgroup. In the case when $G$ is finitely presented, the accessibility theorem of Dunwoody further guarantees that $G$ decomposes as a finite graph of groups with finite edge-groups and one-ended vertex groups.

In the context of this work it is important to stress the topological interpretation of a group being one-ended. It is not hard to show that $e(X)$ is a quasi-isometry invariant. Therefore, $e(G)=e(X)$ for all geodesic metric spaces admitting a geometric group action by $G$. If now $\hat X$ is \emph{any} compactification of $X$ containing $X$ as an open subspace -- call such compactifications \emph{good}, -- then the boundary $\bd X=\hat X\minus X$ is necessarily connected, being the intersection of a descending chain of connected compact subset of $\hat X$. Thus, in order for $G$ to have more than one end (and, consequently, to split over a finite subgroup) it is necessary and sufficient that some geodesic space admitting a geometric action by $G$ have a good compactification with disconnected boundary.

To summarize, in order to employ topological methods for the study of relations between asymptotic/coarse properties of $G$ and its various splittings it will be enough, in view of Dunwoody accessibility, to consider one-ended finitely presented groups $G$.

At this point it is noteworthy to emphasize that the various techniques employed by different authors to provide ``geometric proofs'' of Stallings' theorem, as well as other splitting results, have all utilized the same result by Bass and Serre: a group $G$ splits if and only if it acts on a directed simplicial tree $T$ without a global fixed vertex and without a global fixed end.

The proof of this theorem provides a tool allowing us, once we are able to control the edge-stabilizers of $T$, to have precise information regarding the nature of the splittings arising from $T$.\\

Unfortunately, the approaches that worked for the relatively-hyperbolic case (\cite{[Bo1],[Gu]}) turned out to be hard to apply to the CAT(0) setting, because the techniques used depend strongly on the local connectivity of the boundaries involved, as well as on special dynamic properties that the action of a relatively-hyperbolic group on its canonical boundary has (for example, one question -- still open --, asks when is the action of a group on a CAT(0) boundary minimal). The few known examples of CAT(0) groups whose boundaries seem to portray information regarding known splittings of these groups are also known to have non-locally-connected multiple distinct boundaries (see \cite{[MiRu1],[MiRu2],[CroKle]}), and the methods that were originally developed for (relatively) hyperbolic groups fail. Therefore, it seemed a more immediate goal to try and pin-down those connectivity properties of a CAT(0) group which turn up in \emph{all} of its boundaries. In order to do so, we turned back to analyzing the end structure of $G$.

The main problem realizing the same approach for splittings of one-ended groups over infinite subgroups has been the absence of a good combinatorial object whose structure expresses what one would like to consider as the ``relative end structure'' of the pair $(G,H)$.

While in the case when $H$ is the trivial subgroup of $G$ and the pair $(G,H)$ is multi-ended one naturally expects (from the point of view of coarse geometry) the ``end structure'' of the pair $(G,H)$ to be the same as that of a tree, when dealing with infinite subgroups $H$ this is not what one gets. Instead, in his thesis (\cite{[Sa]}) Sageev has shown that what one usually gets is an action on a higher-dimesional analog of a tree -- namely, a \emph{cubing}, or \emph{non-positively curved cube complex}. Still, the action one obtains is non-trivial in a sense similar to the non-triviality condition of Bass and Serre, so that one may hope to find an equivariant retraction of this cube complex onto a tree in order to achieve a splitting. The main problems with this approach are:
\begin{itemize}
    \item[-] Given the subgroup $H$, the ``Sageev cubing'' is not uniquely determined by $H$; even slight modifications in the choices involved in Sageev's construction may change the dimension of the complex drastically;
    \item[-] Unlike the case of finite subgroups, it is not at all clear whether two different Sageev cubings of the same pair $(G,H)$ have the same ``trace'' on the boundary of a given compactification of $X$.
    \item[-] Unlike the tree constructed for the proof of Stallings theorem, or the end compactification of $\Gamma$, there is no notion of boundary for cubings that will be coarse enough to enable comparison of coarse properties of Sageev cubings with coarse properties of $\Gamma$ (or of any geodesic space $X$ admitting a geometric action by $G$).
    \item[-] There is no known bound on the dimension of the Sageev cubings.
\end{itemize}

Of these, the main obstacles seemed to be the second and third. Thus, the first part of this work deals with constructing a combinatorial boundary for cubings and with its properties. This involves deepening the analysis of Roller's duality between cubings (in their incarnation as discrete median algebras) and $\omega$-dimensional poc-sets (surveyed in section \ref{section:prelim}) by investigating hierarchical relations among the distinct components of the double-dual of a cubing (section \ref{subsection:ae of UFs}). In honor of Roller's pioneering work on this duality, we call the resulting structure by the name \emph{the Roller boundary of a cubing} (defined in \ref{The Roller boundary defined}).\\

Because of the abundance of open questions regarding the Sageev cubings of a group pair, the main test-case for our construction in this thesis was chosen to be less problematic. We assume our group $G$ acts geometrically on a CAT(0) space $X$, with $X$ containing a fixed $G$-invariant \emph{halfspace system}. The idea had been -- before any attempts at coarse geometry are made -- to apply the Roller boundary to better tailored objects, like hyperplanes in Euclidean or Hyperbolic space, or walls in the Davis-Moussong complex of a Coxeter group. We have borrowed the notion of a wall and of a halfspace from these geometric examples, obtaining the following definitions: a \emph{halfspace} in a CAT(0) space $X$ is a convex open subset $h\subset X$ such that $h^\ast=X\minus\cl{h}$ is also convex; the \emph{wall} $W(h)$ of a halfspace $h$ is then defined to be $W(h)=\cl{h}\cap\cl{h^\ast}$; a \emph{halfspace system} is then a set $\HH$ of halfspaces which is closed under the operation of replacing $h$ by $h^\ast$, and satisfying some obvious discreteness conditions, like, for example, that of any point in $X$ having a neighbourhood meeting only finitely-many walls of $\HH$.\\

It is a crucial point in this work that one can use Sageev's construction to show that a halfspace system in CAT(0) space naturally defines a cubing $\hsm{C}(\HH)$. If the halfspace system $\HH$ is $G$-invariant, then many well-known examples (see \cite{[Sa]} for a general discussion, for Coxeter groups see \cite{[NibRee1],[Wil]}, for small cancellation groups see \cite{[Wise]}) demonstrate the importance of knowing when is the action of $G$ on this cubing co-compact. In this work, we answer this question in terms of properties of a canonical map $\rho_{\HH}$ we have defined from the CAT(0) boundary of $X$ into the Roller boundary of $\hsm{C}(\HH)$, thus showing that this question is of a coarse nature (see theorem \ref{thm:co-compactness criterion}).\\

Apart from this main result, we show several other uses for the map $\rho$. In particular, we show that $\rho$ induces a stratification of $\bd X$ with respect to the cone topology (proposition \ref{prop:closure formula 0}), modeled on the hierarchical structure of the Roller boundary $\Re\HH=\Re C(\HH)$, and that the same stratification is inherited by convex subspaces of $X$ (proposition \ref{prop:the restriction equality}).

Under additional geometric assumptions on $\HH$, which are satisfied by CAT(0) cube complexes as well as by Davis-Moussong complexes of Coxeter groups of finite rank, we provide a formula calculating the closure of a stratum (theorem \ref{thm:closure formula 1}), prove the injectivity of $\rho$ on subsets of $\bd X$ which are $\pi$-discrete in the angular metric (corollary \ref{cor:rho separates pi-discrete sets}), and show how the decomposition into strata computes the Tits path-components of $\bd X$ (proposition \ref{prop:safe paths} and theorem \ref{thm:Tits components computed}). In particular, we show how our tools work for several known examples of groups and their CAT(0) boundaries, such as the Croke-Kleiner example.\\

\textbf{Acknowledgements.} some of the results herein were obtained, though in a slightly weaker form, as parts of the author's Ph.D. thesis. The generous support of the Technion IIT is gratefully acknowledged. Many thanks go to my scientific advisors Michah Sageev and Bronislaw Wajnryb, for giving me all the support and guidance I could ever hope for.

\section{Preliminaries: Sageev-Roller duality.}\label{section:prelim}
A cubing is a piecewise-Euclidean simply-connected cell complex
$\hsm{C}$ satisfying the following requirements:
all cells of $\hsm{C}$ are standard Euclidean cubes;
all attaching maps are Euclidean isometries;
no two $k$-faces of the same $d$-cube are attached to each other, for all $0\leq k\leq d$;
all links are (simplicial) flag-complexes (Gromov's so called \emph{link condition} for non-positive curvature).

Given a cubing $\hsm{C}$, it is possible to make $\hsm{C}$ a length-metric space using the identification of each $d$-dube with the unit Euclidean $d$-cube. By a theorem of
Bridson (\cite{[BH]}, theorems I-7.19,I-7.50) and the so-called Cartan-Hadamard theorem (\cite{[BH]} theorem II-4.1, combined with II-5.20), if $\hsm{C}$ is finite-dimensional, or if $\hsm{C}$ is locally finite, then the induced path pseudo-metric on $\hsm{C}$ is a complete CAT(0) metric and every cube of $\hsm{C}$ is isometrically embedded in $\hsm{C}$.

\subsection{The halfspace structure of a cubing.}\label{walls in cubing} In a cubing $\hsm{C}$, the notion of a wall/halfspace arises naturally in combinatorial form. Two edges in a cubing are said to be parallel, if there exists a $2$-cube containing them as opposite edges; one extends this notion of parallelism over the $1$-skeleton of $\hsm{C}$ by taking its transitive closure; next, given a $d$-dimensional cube $Q$ of $\hsm{C}$, one may divide its $1$-skeleton into $d$ disjoint classes of edges which are parallel in $Q$, and define a \emph{midplane} of $Q$ with respect to a class $E$ (computed in $Q$ -- not in $\hsm{C}$) to be the convex hull in $Q$ of the set of midpoints of edges in $E$. Note that if $M$ is a midplane of a cube $Q$, and $Q'$ is a face of $Q$, then $M\cap Q'$ is a midplane of $Q'$ whenever this intersection is non-empty. If $Q_1,Q_2$ are adjacent cubes, and $Q'$ is their maximal common face, then the midplanes $M_i$ of $Q_i$ are said to be compatible if and only if the sets $M_i\cap Q'$ are non-empty and equal. Once again extending the compatibility relation transitively, we say that the union of a compatibility class of midplanes in $\hsm{C}$ is \emph{a wall}. The following is proved in Sageev's thesis:

\begin{thm}[\cite{[Sa]}, theorems 4.10, 4.11] Suppose $\hsm{C}$ is a cubing and $W$ is a wall of $\hsm{C}$. Then:
\begin{enumerate}
    \item $W$ is itself a cubing;
    \item $W$ does not self-intersect -- that is: the intersection of $W$ with any cube of $\hsm{C}$ either equals the empty set or a unique midplane of that cube.
    \item $W$ separates $\hsm{C}$ into the union of precisely two connected componets, whose common boundary in $\hsm{C}$ equals $W$ -- these components are called \emph{the halfspaces determined by $W$}.\ep
\end{enumerate}
\end{thm}
\textbf{Note: } One may also extend the discussion in the proof of theorem 4.13 of \cite{[Sa]} or use Roller's duality results in order to show that the halfspaces of a cubing are convex.\\

\subsection{Duality of poc-sets and cubings.} It has been Sageev's major discovery that a cubing may be fully reconstructed from the nesting patterns of its halfspaces. We present the construction in the form introduced by Roller in \cite{[Rol]}, as an application of a duality theory for median algebras.

\begin{defn}[poc-set, nesting, transversality]\label{defn:poc-set terminology} A poc-set $(H,\leq,\ast)$ is a partially-ordered set $(H,\leq)$ with a minimum element $0$ and an order-reversing involution $h\mapsto h^\ast$ satisfying the requirement that for all $h\in
H$, if $h\leq h^\ast$ then $h=0$. 
\begin{itemize}
    \item[-] the elements $0,0^\ast$ are said to be the \emph{trivial} elements of $H$, while all other elements of $H$ are \emph{proper}.
    \item[-] the poc-set $(H,\leq,\ast)$ is said to be \emph{discrete}, if, for every pair of proper elements $a,b\in H$, the interval $[a,b]=\{h\in H\,|\,a\leq h\leq b\}$ is finite.
    \item[-] two elements $h,k\in H$ are said to be nested (resp. transverse), -- denoted here with $h\| k$ (resp. $h\pitchfork k$) -- if one (resp. none) of the relations $h\leq k,\,h^\ast\leq k,\,h\leq k^\ast,\,h^\ast\leq k^\ast$ holds. A subset $S\subseteq H$ is nested (resp. transverse) if all its elements are pairwise nested (resp. transverse).
    \item[-] a poc-set $(H,\leq,\ast)$ is said to have dimension at most $d$ if every transverse subset of $H$ has at most $d$ elements. $H$ is said to be $\omega$-dimensional, if it contains no infinite transverse subset.
\end{itemize}
\end{defn}
A basic observation by Sageev is that the halfspace system $\HH$ of a cubing $\hsm{C}$ is a discrete poc-set with respect to inclusion and the complementation operator defined by $h^\ast=\hsm{C}\minus\cl{h}$.

\subsubsection{Ultrafilters.} Suppose now we are given a discrete poc-set $H$, and we wish to construct a cubing $\hsm{C}$ whose natural halfspace system $\HH$ is poc-isomorphic to $H$ (with the obvious definition of a poc-morphism as a $(\ast)$-equivariant morphism of posets). One constructs a dual space for $H$:
\begin{defn}\label{defn:ultrafilter} Suppose $(H,\leq,\ast)$ is poc-set. An ultrafilter $\alpha$ on $H$ is a subset of $H$ satisfying:
\begin{description}
    \item[$\mathbf{(UF1)}$] for all $h\in H$, either $h\in\alpha$ or $h^\ast\in\alpha$, but not both;
    \item[$\mathbf{(UF2)}$] for all $h,k\in\alpha$, the relation $h\leq k^\ast$ is prohibited.
\end{description}
The space of all ultrafilters on $H$ will be denoted by $H^\circ$.

A collection $\alpha\subset H$ satisfying $\mathbf{(UF2)}$
is called a \emph{filter base}.
\end{defn}
\begin{remark} Zorn's lemma implies that any filter base is contained in an ultrafilter.
\end{remark}
\begin{remark} $H^\circ$ inherits a topology from $2^H$ (the Tychonoff topology). In fact, $H^\circ$ is a Stone space with respect to this topology.
\end{remark}
One motivation for considering ultrafilters on a poc-set is that when $H$ is the natural halfspace system of a cubing $\hsm{C}$, the vertices of $\hsm{C}$ may be mapped naturally into $H^\circ$ by sending every vertex $v$ to the ultrafilter $\pi_v$ of all halfspaces containing $v$. It is clear that ultrafilters arising in this way satisfy a \emph{descending chain condition}:
\begin{defn}\label{defn:prinicipal UF}
An ultrafilter $\alpha$ on a poc-set is \emph{principal}, if it satisfies the \emph{descending chain condition} (DCC): if $(h_n)_{n=1}^\infty$
is a descending chain of elements in $\alpha$ in the sense that
$h_{n+1}\leq h_n$ for all natural $n$, then $h_{n+1}=h_n$ for all
but finitely-many values of $n$.
\end{defn}
Roller considered a stronger condition for the purpose of distinguishing the ultrafilters of the form $\pi_v$:
\begin{defn}\label{defn:well-founded UF} An ultrafilter $\alpha$ on a poc-set $H$ is {\emph well-founded}, if, for every $a\in\alpha$, $\alpha$ contains only finitely many elements $h$ satisfying $h\leq a$.
\end{defn}

\subsubsection{Almost-equality and (re)constructing cubings.} Since $H$ is discrete, it makes sense to consider the following metric on $H^\circ$ (infinite values are allowed):
\begin{equation}
	\Delta(\alpha,\beta)=\frac{1}{2}\big|\alpha\vartriangle\beta\big|=\left|\alpha\minus\beta\right|=\left|\alpha\cap\beta^\ast\right|\,,
\end{equation}
When $\Delta(\alpha,\beta)$ is finite, we say that $\alpha,\beta$ are {\it almost equal}. For example, if $H$ is the halfspace system of a cubing, then the ultrafilters $\pi_v$ described above all lie in the same almost-equality class.

The metric $\Delta$ allows one to construct a graph $\Gamma=\Gamma(H)$ with vertex set $H^\circ$ and with $\alpha,\beta\in H^\circ$ joined by an edge iff $\Delta(\alpha,\beta)=1$. It is clear that the components of $\Gamma$ correspond to almost-equality classes of $H^\circ$. Henceforth, we will use $\Gamma_\Sigma$ to denote the component of an almost-equality class $\Sigma$ in the graph $\Gamma$. A more detailed study of these graphs will be done in the next section, while here we focus on how they are used for constructing cubings.\\

Sageev (in a special case -- cite{[Sa]}) and Roller \cite{[Rol]} showed that each of the graphs $\Gamma_\Sigma$ is the $1$-skeleton of a cubing:
\begin{description}
	\item[$1$-skeleton.] An important element of their observation was that the restriction of $\Delta$ to a class $\Sigma$ coincides with the combinatorial metric $\Sigma$ induced from $\Gamma_\Sigma$. This facilitated the following inductive construction of a cubing $C_\Sigma$ whose $1$-skeleton is $C_\Sigma^1=\Gamma_\Sigma$.
	\item[$d$-skeleton.] Given $d\geq 2$ and a piecewise-Euclidean cubical $(d-1)$-complex $C_\Sigma^{d-1}$ (consisting of unit cubes), one glues a unique isometric copy of $[0,1]^d$ onto every occurrence of $\partial[0,1]^d$ in $C_\Sigma^{d-1}$. The resulting complex is $C_\Sigma^d$.
\end{description}
Thus, if one wants to reconstruct a cubing $C$ from its halfspace system, one `only' needs to select the right almost-equality class. Sageev and Roller have shown that the almost-equality class -- call it $\Pi$ -- of the ultrafilters $\pi_v$ defined above is the right choice, in the sense that the vertex map $v\mapsto\pi_v$ induces an isomorphism of cubings $C\to C_{\Pi}$.\\

In the general situation (when $H$ is provided abstractly) it is harder to select the right almost-equality class. And what is the meaning of `right' anyway? In \cite{[Gu-loc-fin]} the author shows that both the sets of principal ultrafilters and of well-founded ultrafilters form unions of almost-equality classes. It can also be shown (see \cite{[Rol]}, proposition 9.4) that every well-founded almost-equality class is Tychonoff-dense in $H^\circ$, and one easily sees that $H$ can be recovered from the corresponding cubing (although the different cubing arising in this way may not be isomorphic). Thus, knowing $H$ is the halfspace system of a cubing $C$, provides a way for associating a canonical dense almost-equality class of $H^\circ$ with $H$ the cubing $C$.

In general though, one needs to know more about the poc set $H$ in order to select a `nice' canonical class. In \cite{[Gu-loc-fin]} it is shown that a canonical class of principal ultrafilters exists in either of the following cases:
\begin{description}
	\item[$H$ has a principal ultrafilter of finite degree in $\Gamma$. ] Corollary \ref{cor:misc1} shows that such a vertex of $\Gamma$ is a well-founded ultrafilter, and a computation shows that {\it in this case} there can be only one well-founded class.
	\item[$H$ is $\omega$-dimensional. ] In this case there is only one almost-equality class of principal ultrafilters, which is then also well-founded (corollary \ref{only one principal class}).
\end{description}
In either case, we denote the obtained canonical class by $\Pi$, and call it the {\it principal class} of $H^\circ$. In the geometric application discussed in this paper the first condition will be fulfilled automatically, while for some of the applications we will have to invoke the additional assumption of $\omega$-dimensionality on $H$.

\subsubsection{Generalities about the structure of the graph $\Gamma(H)$.}
Our analysis of cubings will require some technical details. The following observations relating to the local structure of the graph $\Gamma$ will be used frequently, and appear already in \cite{[Sa]}.\\

Given an edge of $\Gamma$ with endpoints $\alpha,\beta\in H^\circ$, it is clear that $\beta$ has the form
\begin{equation}
	[\alpha]_a=(\alpha\minus\{a\})\cup\{a^\ast\}
\end{equation}
for some $a\in\alpha$. It is readily seen that $\beta$ being an ultrafilter implies $a\in\min(\alpha)$: by this we mean that $a$ has to be a minimal element in $\alpha$, where $\alpha$ is considered with the ordering induced from $H$. Thus
\begin{lemma} The set of neighbours of any given $\alpha\in H^\circ$ is parametrized by $\min(\alpha)$.
\end{lemma}
For example, it immediately follows that --
\begin{cor}\label{cor:misc1} If $\alpha\in H^\circ$ is principal and has finite degree then $\alpha$ is well-founded.
\end{cor}

Another case of interest arises when $\alpha$ is the vertex of a $d$-dimensional cube $Q$: in this case, if $a_1,\ldots,a_d$ are the minimal elements of $\alpha$ corresponding to its neighbours in $Q$, then, for every pair of distinct $i,j$ one must have 
\begin{equation}
	\left[[\alpha]_{a_i}\right]_{a_j}=\left[[\alpha]_{a_j}\right]_{a_i}
\end{equation}
which implies $a_i\pitchfork a_j$ for all such pairs. We obtain --
\begin{lemma} For any $\alpha\in H^\circ$, there is a one-to-one correspondence between transverse subsets of $\min(\alpha)$ of size $d$ and the $d$-dimensional cubes $Q$ arising from the Sageev construction and satisfying $\alpha\in Q^0$.
\end{lemma}
In particular, if $H$ is a $d$-dimensional discrete poc-set, then $C_\Sigma$ is at most $d$-dimensional, for every almost-equality class $\Sigma$ of $H^\circ$; if $H$ is discrete and $\omega$-dimensional, then $C_\Sigma$ contains no infinite dimensional cube. The converse is, unfortunately, false.

\subsubsection{Convexity structure and topology on $H^\circ$.}\label{convexity structure on H circ}
Here we provide a brief account on convexity results from \cite{[Rol]} that are used in this paper.\\

Recall that $2^H$ -- the power set of $H$ -- carries a natural \emph{median operation}:
\begin{equation}
    \med(\alpha,\beta,\gamma)=
    (\alpha\cap\beta)\cup
    (\beta\cap\gamma)\cup
    (\alpha\cap\gamma),
\end{equation}
making $2^H$ into a \emph{median algebra}. It is straightforward to verify that $H^\circ$ is closed under this operation. This allows one to define a notion of convexity in $H^\circ$: given $\alpha,\beta\in H^\circ$, \emph{the interval with endpoints $\alpha,\beta$} is defined to be
\begin{equation}
    [\alpha,\beta]=\left\{\gamma\in
    H^\circ\,\big|\,\gamma=\med(\alpha,\beta,\gamma)\right\};
\end{equation}
a subset $\Sigma$ of $H^\circ$ is said to be \emph{convex}, if it contains $[\alpha,\beta]$ for all $\alpha,\beta\in\Sigma$. A convex subset of $H^\circ$ is a {\it halfspace}, if its complement is also convex.\\

When $H$ is a discrete poc-set, the following subsets of $H^\circ$ are \emph{the halfspaces} of $H^\circ$ -- 
\begin{equation}
    S_h=\left\{\alpha\in H^\circ\,\big|\,h\in\alpha\right\},
\end{equation}
where $h$ ranges over the whole of $H$. It is shown in \cite{[Rol]} that a closed subset of $H^\circ$ is convex iff it is the intersection of a family of halfspaces.\\

In any median algebra -- in particular, in $H^\circ$ -- we have the following result:
\begin{thm}[Helly's Theorem (\cite{[Rol]})] Let $C_1,\ldots,C_m$ be pairwise intersecting convex subsets of a median algebra $M$. Then $\bigcap_{i=1}^m C_i$ is non-empty.
\end{thm}
The common application of Helly's theorem in this work will be the following:
\begin{cor}\label{cor:Helly's theorem} Suppose $B\subset H$ is a finite filter-base in $H$, and $\Sigma$ is an almost-equality class of $H^\circ$ such that $\Sigma\cap S_b$ is non-empty for all $b\in B$. Then $\Sigma\cap\bigcap_{b\in B} S_b$ is non-empty.
\end{cor}
The proof of the corollary relies on the following simple fact:
\begin{lemma} Every almost-equality class of $H^\circ$ is a convex subset of $H^\circ$.
\end{lemma}
\proof{} Fix an almost-equality class $\Sigma$. In order to verify our claim, we need to take any $\sigma\in\HH^\circ$, a pair of ultrafilters $\alpha,\beta\in\Sigma$, and check that $\mu=\med(\alpha,\beta,\sigma)$ lies in $\Sigma$. Let us write
\begin{equation}
    \mu=(\alpha\cap\beta)\cup(\alpha\cap\sigma)\cup(\beta\cap\sigma),
\end{equation}
and observe that
\begin{equation}
    \mu\minus\alpha\subseteq (\beta\cap\sigma)\minus\alpha\subseteq\beta\minus\alpha.
\end{equation}
Since the right-hand side is finite, we are done.\ep\\

\proof{of corollary \ref{cor:Helly's theorem}} The proof is straightforward now: the family $\{\Sigma\}\cup\{S_b\}_{b\in B}$ satisfies the assumptions of Helly's theorem, so its total intersection is non-empty.\ep\\

We have already referred to closed subsets of $H^\circ$ and the Tychonoff topology. Let us make this precise. We topologize $H^\circ$ as a subspace of $2^H$ with the product topology. We will refer to this topology as the Tychonoff topology on $H^\circ$. It is easy to see that the family of all subsets of $H^\circ$ of the form
\begin{equation}
	S_h=\left\{\alpha\subseteq H\,\big|\,h\in\alpha\right\}
\end{equation}
is a sub-base for the Tychonoff topology $H^\circ$. Thus, if $A\subseteq H$ is any finite set, then $V(A)=\bigcap_{a\in A}S_a$ is a (possibly empty) basic open set (for $H^\circ$); by Helly's theorem, for $V(A)$ to be non-empty it is necessary and sufficient that no relation of the form $a\leq b^\ast$ holds for $a,b\in A$, so the family of all $V(A)$ with $A$ ranging over all finite filter-bases of $H$ constitutes a basis for the Tychonoff topology on $H^\circ$.\\

The space $H^\circ$ is obviously Hausdorff. Moreover, it is totally disconnected, since the sets $S_h$ are all clopen. To see that $H^\circ$ is compact it is enough to verify it is closed in $2^P$ (we leave this as an exercise to the reader). It follows that
\begin{prop} $H^\circ$ is a Stone space. Moreover, if $\hsm{F}$ is a family of pairwise-intersecting closed convex subsets of $H^\circ$, then $\bigcap\hsm{F}\neq\varnothing$.
\end{prop}

\subsubsection{Functorial properties.}\label{prelim:functorial properties}
Suppose now that $H_1,H_2$ are poc-sets. A morphism of poc-sets is, by definition, a $(\ast)$-equivariant morphism of posets $f:H_1\to H_2$ satisfying $f(0)=0$. It is easy to see that, given a poc-set $H$, the space $H^\circ$ can be naturally identified (using characteristic functions) with the space of poc-morphisms into the $2$-element poc-set. As a result, given $f$ as above, there is a natural pullback map $f^\circ:H_2^\circ\to H_1^\circ$. For our purposes it will be important to state the following facts from \cite{[Rol]}, sections 3 and 4:
\begin{prop} The assignment $H\mapsto H^\circ$, $f\mapsto f^\circ$ is a contravariant functor of the category of poc-sets with poc-morphisms into the category of Stone median algebras with continuous median-preserving maps.
\end{prop}
The pullback in the poc-set category is exact in the following sense:
\begin{prop} Let $f:H_1\to H_2$ be a morphism of poc-sets.
\begin{itemize}
	\item $f^\circ$ is injective iff $f$ is surjective;
	\item $f^\circ$ is surjective iff $f$ is an embedding.
\end{itemize}
\end{prop}

\section{A boundary for $\omega$-dimensional cubings}\label{section:Roller boundary}
Throughout this section, $(H,\leq,\ast)$ is a discrete $\omega$-dimensional poc-set. Our initial goal will be to establish some elementary facts about how the combinatorics of infinite descending chains is related to the structure of a given almost-equality class and its Tychonoff closure. We first study the principal class.

\subsection{Characterizing the principal class.}
Consider the following way of constructing an ultrafilter:
\begin{prop}\label{prop:constructing principal ultrafilters} Let $A$ be a maximal transverse subset of $H$. Then the collection
\begin{equation}
    \pi_A=\left\{h\in H\,\big|\,
        \exists_{a\in A}\;a\leq h
    \right\}\cup
    \left\{h\in H\,\big|\,
        \exists_{a\in A}\;a^\ast<h
    \right\}
\end{equation}
is a principal ultrafilter on $H$. Moreover, every principal ultrafilter is of this form.
\end{prop}
\proof{} First let us show $\pi_A\in H^\circ$. We need to show three things:
\begin{itemize}
    \item For any $h\in H$, either $h\in\pi_A$ or $h^\ast\in\pi_A$.

    Given $h\in H$, if neither $h$ nor $h^\ast$ are
    greater than an element of $A\cup A^\ast$, it means that
    $A\cup\{h\}$ is a transverse set, contradicting the maximality of
    $A$ as such a set.
    
    \item For any $h\in H$, if $h\in\pi_A$ then $h^\ast\notin\pi_A$.

    If \emph{both} $h$ and $h^\ast$ lie in $\pi_A$, then there are
    $a,b\in A\cup A^\ast$ such that $a\leq h$ and $b\leq h^\ast$.
    This implies $a\leq h\leq b^\ast$, and the transversality of $A$ then gives
    $a=b^\ast$, resulting in $a=h$. Thus, both $h=a$ and
    $h^\ast=a^\ast$ lie in $\pi_A$, which is impossible, by the construction of $\pi_A$.
    
    \item If $h,k\in\pi_A$, then $h\not\leq k^\ast$.

    Suppose the contrary, and let $a,b\in A\cup A^\ast$ be such
    that $a\leq h$ and $b\leq k$. Then we have
    \begin{equation}
        a\leq h\leq k^\ast\leq b^\ast,
    \end{equation}
    meaning again that $b^\ast=a$, and hence $h=k^\ast$, finally
    implying $k,k^\ast\in\pi_A$, which we have just shown to be impossible.
\end{itemize}
In order to prove $\pi_A$ is principal, assume $\eta=(h_n)_{n=1}^\infty$ is a descending chain in $\alpha$, and we have to show it is eventually constant. Indeed, for each $n$ find $a_n\in A\cup A^\ast$ such that $a_n\leq h_n$. Since $A\cup A^\ast$ is finite, the $a_n$ may be taken to be constant -- denote them by $a$. Thus $h_n\in[a,h_1]$ for all $n$. Since $H$ is discrete, the $h_n$ stabilize.\\

Conversely, suppose $\alpha$ is a principal ultrafilter, and consider the set $\min(\alpha)$ of all minimal elements of $\alpha$. Since $\alpha$ is principal, it is clear that every $h\in\alpha$ has some $a\in\min(\alpha)$ satisfying $a\leq h$. In particular, $\min(\alpha)$ is non-empty.

We will now show $\min(\alpha)$ contains a maximal transverse set $A$ of $H$; in particular, $\alpha$ contains $\pi_A$, which is an ultrafilter. Since no two distinct ultrafilters contain each other, we will conclude $\pi_A=\alpha$, completing the proof.

We now construct the required set $A$ by induction. Start from any element $a_1\in\min(\alpha)$, and let $A_1=\{a_1\}$. If $A_1$ is a maximal transverse subset of $H$, then we are done -- else, proceed by induction as follows: suppose $A_i\subset\min(\alpha)$ is a given transverse subset which is not maximal in $H$, then there is an element $h$ of $H$ with $h\pitchfork A_i$, and then there is an element $a\in\min(\alpha)$ with either $a\leq h$ or $a\leq h^\ast$; in any case, $a\pitchfork A_i$, and we set $A_{i+1}=A_i\cup\{a\}$. Since $H$ contains no infinite strictly-ascending chain of transverse subsets, this process must stop, producing a subset $A$ of $\min(\alpha)$ which is a maximal transverse subset of $H$, as desired.\ep

\begin{cor}\label{cor:principality by waist} $\alpha\in H^\circ$ is principal iff $\min(\alpha)$ contains a maximal transverse subset $A$ of $H$, in which case $\alpha=\pi_A$.
\end{cor}
An important consequence of this is:
\begin{cor}\label{only one principal class}
The principal ultrafilters on $H$ form an almost-equality class in $H^\circ$.
\end{cor}
\proof{} Suppose $\alpha=\pi_A$ and $\beta=\pi_B$ are principal ultrafilters, and consider $h\in\alpha\minus\beta$: we have $a\leq h$ for some $a\in A\cup A^\ast$, but $h\notin\beta$; since $B$ is a maximal transverse set, $B\cup\{h\}$ is not transverse, there is a $b\in B$ such that the pair $\{b,h\}$ is nested. The assumption $h\notin\pi_B$ then forces either $h<b$ or $h\leq b^\ast$, so that $h$ belongs to $[a,b]\cup[a,b^\ast]$. Thus, we have shown that $\beta\minus\alpha$ lies in the union of all intervals of the form $[a,b]$ with $a\in A\cup A^\ast$ and $b\in B\cup B^\ast$ whenever those are defined; since $H$ is discrete and both $A$ and $B$ are finite, we are done.\ep\\

Henceforth, we shall denote the set of principal ultrafilters on $H$ by $\Pi(H)$, or simply by $\Pi$, when $H$ is understood from the context.

\begin{cor}\label{cor:topological characterization of principality} $\Pi$ is dense in $H^\circ$. Moreover, it is the only almost-equality class of $H^\circ$ with this property. \end{cor}
\proof{} Let us prove $\Pi$ is dense. Let $V(A)$ be an arbitrary basic open set in $H^\circ$, with $A$ a finite filter base in $H$. Then the corollary (\ref{cor:Helly's theorem}) to Helly's theorem simply states that $V(A)\cap\Pi$ is non-empty, provided we can show that $\Pi\cap S_a$ is non-empty for all $a\in A$. In fact, $\Pi\cap S_a\neq\varnothing$ for any proper $a\in H$: given $a$, let $A$ be any maximal transverse subset containing $a$. Then $\pi_A$ lies in $\Pi\cap S_a$, by proposition \ref{prop:constructing principal ultrafilters}. 

We defer the proof of the second part until after we have studied some combinatorial properties of almost-equality classes.\ep

\subsection{Almost-equality of ultrafilters}\label{subsection:ae of UFs}

Given an almost-equality class $\Sigma$, recall that the metric $\Delta$ on $H^\circ$ restricts on $\Sigma$ to the combinatorial metric of the graph $\Gamma_\Sigma$ (this result, initially due to Sageev in a special case, was proved in this generality, independently in \cite{[Rol]}, \cite{[Nic]} and \cite{[ChaNib]} with slight variations of context).

Thus, a vertex path $(\sigma_0,\ldots,\sigma_d)$ from $\alpha\in\Sigma$ to $\beta\in\Sigma$ must satisfy $d\geq\Delta(\alpha,\beta)$, and may be chosen so that equality holds. Let us recall that for each $i\in\{1,\ldots,d\}$ there are $h_i\in\min(\sigma_i)$ such that $\sigma_{i+1}=\left[\sigma_i\right]_{h_i}$. It will be convenient to consider such a vertex path as a chain (of length $d$) of transformations of $\alpha$ into $\beta$, performed in stages. For this reason, we shall also refer to the (directed) edges of $\Gamma_\Sigma$ as {\it elementary moves}.

\subsubsection{Parallel/Transverse decomposition.}
Let us fix an almost-equality class $\Sigma$ of $H^\circ$. Obviously, some elements of $H$ have nothing to do with separation in $\Sigma$: 
\begin{defn} Suppose $h\in H$ and $\Sigma$ is an almost-equality class in $H^\circ$. We say \emph{$h$ is transverse to $\Sigma$} (denoted by $h\pitchfork\Sigma$), if both $S_h\cap\Sigma$ and $S_{h^\ast}\cap\Sigma$ are non-empty. Otherwise we shall say that $h$ is \emph{parallel to $\Sigma$}. We shall keep the following notation:
\begin{eqnarray}
    T(\Sigma)&=&\left\{h\in
    H\left|h\pitchfork\Sigma\right.\right\}\cup\{0,0^\ast\}.
\end{eqnarray}
\end{defn}
With the induced ordering, $T(\Sigma)$ is a discrete $\omega$-dimensional poc set in its own right. Its less abstract counterpart is the poc-set
\begin{equation}
	H_\Sigma=\left\{
		S_h\cap\Sigma\,\big|\,
		h\in H
	\right\}
\end{equation}
viewed as a sub poc-set of the power set of $\Sigma$. Note that $H_\Sigma$ is the image of $H$ under the \emph{restriction map}
\begin{equation}
	\left\{\begin{array}{rcl}
		r_\Sigma:H  &\to&       H_\Sigma\\
		h  &\mapsto&   S_h\cap\Sigma\,,
	\end{array}\right.
\end{equation}
mapping $T(\Sigma)$ onto $H_\Sigma$, and sending elements parallel to $\Sigma$ to the trivial elements of $H_\Sigma$. Obviously, $r_\Sigma$ is a morphism of poc-sets.
\begin{lemma}\label{lemma:transverse is minimal} Assume $a\pitchfork\Sigma$. Then there exists an ultrafilter $\gamma\in S_a\cap\Sigma$ such that $a\in\min(\gamma)$. Conversely, if $\gamma\in\Sigma$ and $a\in\min(\gamma)$, then $a\pitchfork\Sigma$.
\end{lemma}
\proof{} Let $\gamma\in S_a\cap\Sigma$ and $\gamma^\ast\in S_{a^\ast}\cap\Sigma$. Then there is a sequence of $\Sigma$-edges from $\gamma$ to $\gamma^\ast$ -- in particular, $a$ is contained in the minimal set of some element of this sequence.

For the converse, simply consider $[\gamma]_a$: this is an element of $\Sigma$ contained in $S_{a^\ast}$. Thus, neither $S_a\cap\Sigma$ nor $S_{a^\ast}\cap\Sigma$ is empty.\ep

\begin{cor}\label{cor:descending chain is parallel} Suppose $\sigma\in\Sigma$ contains an infinite descending chain $c=(c_n)_{n=1}^\infty$. Then $\Sigma\subseteq S_{c_n}$ for all $n$.
\end{cor}
\proof{} Fix $n$, let $c=c_n$, and suppose $c\pitchfork\Sigma$. Lemma \ref{lemma:transverse is minimal} then implies there exists $\xi\in\Sigma$ containing $c$ in its minimal set. In particular, $\sigma\minus\xi$ contains the infinite set $\{c_m\,|\,m>n\}$, contradicting the fact that $\sigma$ and $\xi$ are almost-equal. Thus, $c_n$ is parallel to $\Sigma$, and being an element of $\sigma\in\Sigma$, necessitates $\Sigma\subseteq S_c$, as desired.\ep\\

From this we obtain the topological characterization of the principal class:

\proof{of cor. \ref{cor:topological characterization of principality}} We have shown already that the principal class $\Pi(H)$ is dense in $H^\circ$. Now, suppose that $\Sigma$ is non-principal. Then some $\sigma\in\Sigma$ contains an infinite descending chain, implying there exist proper elements $h\in H$ satisfying $\Sigma\subseteq S_h$. In particular, $\Sigma$ does not intersect the non-empty open subset $S_{h^\ast}$ of $H^\circ$.\ep

\begin{lemma}\label{lemma:restriction to an ae-class} Suppose $a,b\pitchfork\Sigma$. Then, $S_a\cap\Sigma\subseteq S_b\cap\Sigma$ holds if and only if $a\leq b$. Thus, $r_\Sigma$ restricts to an isomorphism between the poc-sets $T(\Sigma)$ and $H_\Sigma$. In particular, $H_\Sigma$ is a discrete $\omega$-dimensional poc-set.
\end{lemma}
\proof{} One direction is trivial. For the other, we assume both $S_a\cap\Sigma\subseteq S_b\cap\Sigma$ and $a\not\leq b$ hold, and seek to derive a contradiction. 

Among all pairs $(\alpha,\beta)\in\Sigma\times\Sigma$ satisfying the condition
\begin{equation}
    (\ast)\quad
    a\in\min(\alpha)\;\wedge\;b^\ast\in\min(\beta)
\end{equation}
let us fix one for which $n=\Delta(\alpha,\beta)$ is minimal. Note that lemma \ref{lemma:transverse is minimal} guarantees the set of all pairs satisfying $(\ast)$ is non-empty.\\

Let $(\sigma_0=\alpha,\ldots,\sigma_n=\beta)$ be a chain (of elementary moves) of length $n$. We define
\begin{eqnarray}
    i&=&\max\left\{t\in\{0,\ldots,n\}\,\big|\,\forall_{s\leq t}\;a\in\sigma_s\right\}\\
    j&=&\min\left\{t\in\{0,\ldots,n\}\,\big|\,\forall_{s\geq t}\;b^\ast\in\sigma_t\right\}.
\end{eqnarray}
Since $a\in\sigma_i$, our initial assumption implies $b\in\sigma_i$, so that $i<n$. Note now that, by the definition of $i$, $a^\ast\in\sigma_{i+1}$, so that $a\in\min(\sigma_i)$. Thus, if $i$ were greater than zero, then the chain $(\sigma_i,\ldots,\sigma_n)$ satisfies $(\ast)$ -- a contradiction to $n$ being the minimal length of such a chain. Similarly, one proves $j=n$, and we have
\begin{eqnarray}
    t=1,\ldots,n    &\THEN& a^\ast\in\sigma_t\\
    t=0,\ldots,n-1  &\THEN& b\in\sigma_t.
\end{eqnarray}
Note that $a^\ast\in\min(\sigma_1)$ and $b\in\min(\sigma_{n-1})$.

Assume $a^\ast\in\min\sigma(t)$, and let us show $a^\ast\in\min(\sigma_{t+1})$. By definition of an elementary move, there is an element $c\neq a^\ast$, $c\in\min(\sigma_t)$ such that $\sigma_{t+1}=[\sigma_t]_c$. If $a^\ast$ is not minimal in $\sigma_{t+1}$, then it means $c^\ast<a^\ast$, which implies $c<a$. But then $c\in\sigma_t$ implies $a\in\sigma_t$, which is impossible, as $a^\ast\in\sigma_t$. 

Thus, we have shown $a^\ast\in\min(\sigma_t)$ for all $t=1,\ldots,n$. In particular, the ultrafilter $\sigma=\sigma_{n-1}$ contains both $a^\ast$ and $b$ as minimal elements, and $a^\ast\not\leq b^\ast$. 

Recall now that $a,b\in\alpha$, so that $a\not\leq b^\ast$. Further, since $a^\ast,b^\ast\in\beta$, we also have $a^\ast\not\leq b$. Adding to this the assumption that $a\not\leq b$, we conclude $a\pitchfork b$ holds. Thus, the ultrafilter $\mu=[\sigma]_{a^\ast,b}$ is an element of $\Sigma$ containing both $a$ and $b^\ast$. Now, applying our assumption that $\Sigma\cap S_a\subseteq \Sigma\cap S_b$ we obtain 
\begin{equation}
    a\in\mu\IFF \mu\in S_a\THEN \mu\in S_b,
\end{equation}
contradicting $\mu\in S_{b^\ast}$.\ep

\subsubsection{Hierarchy in $H^\circ$.}
\paragraph{restriction to a class, corestriction, projection.}
Our aim now is to use the fact that $r_\Sigma$ restricts to an isomorphism of $T(\Sigma)$ with $H_\Sigma$ for computing the Tychonoff closure of $\Sigma$ in $H^\circ$: intuitively it is expected that $\Sigma$ arises as the principal class for each of these poc-sets, and we want to conclude that the structure of the closure of $\Sigma$ in $H^\circ$ depends only on the structure of $\Sigma$.\\

Let us consider the sequence of poc-set morphisms
\begin{equation}
   \begin{CD}T(\Sigma) @>inc.>> H @>r_\Sigma>> H_\Sigma\end{CD}\;.
\end{equation}
We obtain three continuous maps: a map $\iota:H^\circ\to
T(\Sigma)^\circ$ dual to the inclusion, a ``corestriction'' map
$cor^H_\Sigma:H_\Sigma^\circ\to H^\circ$ dual to $r_\Sigma$, and
-- by lemma \ref{lemma:restriction to an ae-class} -- the
homeomorphism $t_\Sigma$ induced (as a push forward) by the
composition of $r_\Sigma$ and the inclusion map. The map
$t_\Sigma$ enables the construction of a ``restriction map''
\begin{equation}
    res^H_\Sigma=t_\Sigma\circ\iota:H^\circ\to H_\Sigma^\circ,
\end{equation}
Explicitly:
\begin{equation}
\begin{array}{rclc}
    res^H_\Sigma(\sigma)&=&
        \left\{S_h\cap\Sigma\,\big|\,h\in\sigma\right\}&
        \left(\sigma\in H^\circ\right)\,,\\
    cor^H_\Sigma(\bar{\sigma})&=&
        \left\{h\in H\,\big|\,S_h\cap\Sigma\in\bar{\sigma}\right\}&
        \left(\bar{\sigma}\in H_\Sigma^\circ\right)\,.
\end{array}
\end{equation}

\begin{prop}\label{lemma:corestriction identities}
Suppose $\Sigma$ is an almost-equality class in $H^\circ$. Then the following hold:
\begin{enumerate}
    \item $res^H_\Sigma\circ cor^H_\Sigma=id_{H_\Sigma}$. In particular, $cor^H_\Sigma$ is injective.
    \item The map $pr_\Sigma:H^\circ\to H^\circ$ defined by $pr_\Sigma=cor^H_\Sigma\circ res^H_\Sigma$ is a retraction, pointwise fixing $\cl{\Sigma}$.
    \item $cor^H_\Sigma$ is a homeomorphism of $H_\Sigma^\circ$ onto $\cl{\Sigma}$.
    \item $cor^H_\Sigma$ sends almost-equality classes to almost-equality classes. In particular, the principal class of $H_\Sigma$ is mapped bijectively onto $\Sigma$.
\end{enumerate}
\end{prop}
\proof{} In the process of this proof, let us suppress the $H$ and $\Sigma$ indices in the notation for the maps in concern, referring to them as $res$, $cor$ and $pr$, respectively.

Claim (1.) is immediate from the explicit forms given above, and implies $pr_\Sigma$ is, indeed, a retraction.\\

By continuity, to prove (2.) it is enough to verify that $pr$ fixes $\Sigma$ pointwise. For $\sigma\in\Sigma$ we compute:
\begin{eqnarray}
    a\in \underbrace{(cor\circ res)(\sigma)}_{\triangleq\sigma'}
        &\IFF& S_a\cap\Sigma\in r_\Sigma(\sigma).
\end{eqnarray}
Now, there are two possibilities: if $a\pitchfork\Sigma$, then lemma \ref{lemma:restriction to an ae-class} tells us that $a\in\Sigma$; if $a$ is parallel to $\Sigma$, then its lying in the image of $r_\Sigma$ means $\Sigma\subseteq S_a$, again implying $a\in\Sigma$. Thus, we have found out that $\sigma'$ is contained in $\sigma$; being an element of $H^\circ$, it is equal to $\sigma$, as desired.\\

For (3.), note that $cor$ is a bijective continuous map onto its image, with $H^\circ$ Hausdorff and $H_\Sigma^\circ$ compact. Thus, it is a homeomorphism onto its image. Therefore, we now have to show $Im(pr)=Im(cor)$ is actually equal to the closure of $\Sigma$, implying $cor$ is a homeomorphism of $H_\Sigma^\circ$ onto $\cl{\Sigma}$. Since $\Pi(H_\Sigma)$ is dense in $H_\Sigma^\circ$ (see cor.\ref{cor:topological characterization of principality}) and $cor$ is continuous, it is enough to show $cor$ maps the principal class $\Pi(H_\Sigma)$ into $\Sigma$. Since $pr=cor\circ res$, this is equivalent to proving $res$ maps $\Sigma$ \emph{onto} $\Pi(H_\Sigma)$. In order to do this, it will be enough to show that $\Sigma$ is mapped into $\Pi(H_\Sigma)$, and then to show that every chain of elementary moves beginning in $res(\Sigma)$ may be lifted through $res$ to a chain of elementary moves in $\Sigma$; of course, it is enough to consider chains of length $1$. Now, for any $\sigma\in\Sigma$, note that corollary \ref{cor:descending chain is parallel} implies $res(\sigma)$ contains no infinite descending chains; thus, $res$ maps $\Sigma$ into the principal class of $H_\Sigma^\circ$, as we expected. To finish our argument, let now $\sigma\in\Sigma$ be any ultrafilter, and consider $a\in H$ such that $\bar{a}=r_\Sigma(a)$ is a minimal element of $res(\sigma)$; in particular, $a\pitchfork\Sigma$. Applying lemma \ref{lemma:restriction to an ae-class} gives us that $a$ is a minimal element of $\sigma$, meaning that both the ultrafilters $[\sigma]_a$ and $[res(\sigma)]_{\bar{a}}$ exist; it is then obvious from the formula we have for $res$ that
\begin{equation}
    (\ast)\qquad res\left([\sigma]_a\right)=[res(\sigma)]_{\bar{a}}.
\end{equation}
We note that this equality is true \emph{whenever the ultrafilters involved are well-defined}, and an analogous equality holds for $cor$; this is important for the proof of (4).

Finally, we prove (4). Assume $\Xi$ is an almost-equality class in $H_\Sigma^\circ$. First of all, (3.) states $cor$ maps $\Xi$ into the closure of $\Sigma$ in $H^\circ$, so that $pr$ is pointwise fixed on $cor(\Xi)$. For every $h\in H$, if $\Sigma\subseteq S_h$, then $\cl{\Sigma}\subseteq S_h$, because $S_h$ is closed. For all $h\in H$ we have
\begin{equation}
    (\ast\ast)\quad\Sigma\subseteq S_h\THEN\Xi\subseteq S_h.
\end{equation}
Applying the same method as in the proof of (3.), in order to prove that $cor(\Xi)$ lies inside an almost-equality class in $H^\circ$, we assume $\xi\in\Xi$ and $a\pitchfork\Sigma$ such that $\bar{a}=S_a\cap\Sigma\in\min(\xi)$ and we need to show that $a$ is minimal in $cor(\xi)$: indeed, if $b\in\xi$ satisfies $b\leq a$, then the minimality of $\bar{a}$ implies $S_b\cap\Sigma=S_a\cap\Sigma$, implying $b\pitchfork\Sigma$; but then lemma \ref{lemma:restriction to an ae-class} shows that $a=b$, and we have shown $a\in\min(\xi)$.

For the reverse inclusion (showing $cor$ maps $\Xi$ \emph{onto} an almost-equality class of $H^\circ$), we need to show that if $a$ is a minimal element in $cor(\xi)$, then $S_a\cap\Sigma$ is a minimal element of $\xi$ (which then enables us to apply $(\ast)$). Now, a minimal element in $cor(\xi)$ is transverse to $\Xi$, which, through $(\ast\ast)$, implies it is also transverse to $\Sigma$. Therefore, if $S_b\cap\Sigma\subseteq S_a\cap\Sigma$ and $S_b\cap\Sigma\in\xi$, then $b\in cor(\xi)$ either satisfies $b\pitchfork\Sigma$ and then $b=a$ (as desired), or $b$ is parallel to $\Sigma$, implying $S_b\cap\Sigma=\varnothing$; but then:
\begin{eqnarray}
    S_b\cap\Sigma=\varnothing
        &\IFF&  \Sigma\subseteq S_{b^\ast}\\
        &\THEN& \Xi\subseteq S_{b^\ast}\\
        &\IFF&  S_b\cap\Xi=\varnothing
\end{eqnarray}
-- a contradiction to $\xi\in S_b\cap\Xi$. Thus, $a$ is a minimal element of $cor(\xi)$, and we are done.\ep\\

For us, an important corollary of the last proposition is the
following
\begin{cor}\label{cor:ordering the Roller boundary} Suppose an almost-equality class $\Xi$ intersects the closure of an almost-equality class $\Sigma$. Then $\Xi$ is contained in the closure of $\Sigma$.
\end{cor}
\proof{} Consider once more the map $cor^H_\Sigma:H_\Sigma^\circ\to H^\circ$, and let $\alpha\in\Xi\cap\cl{\Sigma}$. By the above proposition, $\alpha$ lies in the image of the corestriction map. Moreover, part (4.) of the proposition implies $cor^H_\Sigma$ maps the whole almost-equality class of $\alpha$ into $\cl{\Sigma}$.\ep

\paragraph{Almost-equality and parallelism.}
We have so far considered almost-equality classes disregarding the elements of $H$ parallel to the respective classes. This gave us a kind of ``intrinsic'' approach to analyzing these classes, but could not serve to provide us with ideas regarding the way these classes interact. In the last corollary we have seen that closures of almost-equality classes intersect in almost-equality classes, this resulting in a kind of hierarchy on the set of all classes. In order to study this hierarchy, we need to have a better understanding of the way how, for a given almost-equality class $\Sigma$, the poc-set $T(\Sigma)$ is embedded in $H$. For this purpose, the following lemma is essential:
\begin{lemma}\label{lemma:parallel elements have no minimum} Let $\Sigma$ be an almost-equality class of $H^\circ$, and suppose a proper element $a\in H$ satisfies $\Sigma\subseteq S_a$. Then there exists $b<a$ with the same property.
\end{lemma}
\proof{} Since $a$ is parallel to $\Sigma$, it is not a minimal element of any $\sigma\in\Sigma$. Let us denote by $P$ the set of immediate predecessors of $a$ in $H$, and suppose there is no (proper) $b<a$ such that $\Sigma\subseteq S_b$; then every $p\in P$ is transverse to $\Sigma$. By lemma \ref{lemma:transverse is minimal}, for each such $p$ there exists $\sigma_p\in\Sigma$ such that $p\in\min(\sigma_p)$. 

Let $p,q\in P$ be distinct, and suppose they are nested. Then, since it cannot be that $p<q$ or $q<p$, we must have either $p^\ast<q$ or $p<q^\ast$. Now, if $p^\ast<q$ then the calculation 
\begin{eqnarray}
    p<a &\THEN& a^\ast<p^\ast<q<a,
\end{eqnarray}
yields a contradiction, leaving $p<q^\ast$ as the only possibility. Thus, there exists \emph{no} ultrafilter $\xi\in H^\circ$ containing both $p$ and $q$, which is equivalent to saying that $S_p\cap S_q$ is empty. Now consider $\sigma_p$ and $\sigma_q$: we have
\begin{equation}
    \{p,q^\ast\}\subset\sigma_p\,,\quad\{q,p^\ast\}\subset\sigma_q,
\end{equation}
If $\sigma_0=\sigma_p,\ldots,\sigma_k=\sigma_q$ is a chain of elementary moves -- that is, $\sigma_{t+1}=[\sigma_t]_{a_t}$ for every $t\in\{1,\ldots,k\}$ -- then, since none of the $\sigma_t$ may contain the sets $\{p,q\}$ and $\{p^\ast,q^\ast\}$, we arrive at a contradiction, by the definition of an elementary move. Thus we have proved that $P$ is a transverse set.

By the assumption of $\omega$-dimensionality, $P$ is finite. Also, $P$ is contained in $T(\Sigma)$. Thus, there exists a principal ultrafilter $\bar{\pi}\in\Pi(H_\Sigma)$ containing $r_\Sigma(P)$ in its minimal set, and then $\pi=cor^H_\Sigma(\bar{\pi})$ is an ultrafilter in $\Sigma$ having $P\subseteq\min(\pi)$. Finally, reversing all the elements of $P$ in $\pi$ one-by-one, we obtain an ultrafilter $\pi'\in\Sigma$ with the property that $a\in\min(\pi')$ -- contradiction.\ep

\paragraph{Truncation.}
We define a generalization of elementary moves, which we call \emph{truncation}. The idea behind this terminology is that a non-principal ultrafilter $\sigma$ may be changed into another ultrafilter containing ``a smaller number'' of infinite descending chains than the original by ``cutting-off'' all elements of $\sigma$ below a certain prescribed level.

The truncation operation on $\xi\in\Sigma$ is meant to generalize the operation of reversing the orientation of a minimal element. However, while inverting a minimal element does not change the a.e.-class of an ultrafilter, inverting a whole (descending) chain of halfspaces, say, in an ultrafilter should ``move'' its a.e.-class down with respect to the ordering we had introduced on the set of almost-equality classes. 
\begin{lemma}[truncation]\label{lemma:truncation of ultrafilters}
Let $\xi\in H^\circ$ and $b\in\xi$. Then
\begin{equation}
    [\xi]_b=\xi\minus\{h\in H\,\big|\,h\leq b\}\cup\{h\in H\,\big|\,b^\ast\leq h\}
\end{equation}
is also an element of $H^\circ$. Moreover, $[\xi]_b$ is almost-equal to $\xi$ if and only if $b$ does not belong to an infinite descending chain of $\xi$.
\end{lemma}
\begin{remark}\label{lemma:transitivity of truncation} Suppose $\xi\in H^\circ$, and let $a,b\in\xi$ satisfy $b<a$. Then, $\left[[\xi]_b\right]_a=[\xi]_a$.
\end{remark}
\proof{of the lemma} Let $\Sigma$ be the almost-equality class of $\xi$. We prove the latter assertion first. If $b$ does lie on a descending chain $\{c_n\}_{n=1}^\infty\subset\xi$, then, obviously, $[\xi]_b$ is not in the same almost-equality class as $\xi$, as $c_n\notin[\xi]_b$ for almost all $n$.

For the converse, it is enough to prove $b\pitchfork\Sigma$. Indeed, if $b$ is transverse to $\Sigma$, then there exists $\sigma\in\Sigma$ with $b^\ast\in\min(\sigma)$. But then $\sigma$ is, by the very definition of $[\xi]_b$, almost equal to $[\xi]_b$, so that $[\xi]_b\in\Sigma$, as desired.

Thus, we assume there is no infinite descending chain in $\xi$ through $b$, and we need to show that $b\pitchfork\Sigma$ holds. Now, since $b\in\xi$, we have $\Sigma\cap S_b$ is non-empty (because $\xi$ is in there), and we only need to show that $\Sigma\cap S_{b^\ast}$ is non-empty, too. Assume, on the contrary, the latter intersection is empty; then we have $\Sigma\subset S_b$, and we may use lemma \ref{lemma:parallel elements have no minimum} to inductively construct an infinite descending chain $\{c_n\}_{n=1}^\infty$ such that $c_1=b$ and with the property that $\Sigma\subset S_{c_n}$ for all $n$. In particular, $\xi$ -- being a member of $\Sigma$ -- will contain this chain.\\

Finally, let us show $[\xi]_b$ is an ultrafilter. As before, we break the proof into three parts:
\begin{itemize}
    \item \underline{Suppose $h\in\xi_1$, and show $h^\ast\notin\xi_1$:}

    Assume both $h\in\xi_1$ and $h^\ast\in\xi_1$. For $h$ there are two possibilities with respect to $\xi$:
    \begin{itemize}
        \item If $h\in\xi$, then $h^\ast\notin\xi$, but $h^\ast\in\xi_1$ then implies $b^\ast\leq h^\ast$, giving $h\leq b$, which is impossible for $h\in\xi\cap\xi_1$. 
        \item If $h^\ast\in\xi$, then $h\notin\xi$ together with $h\in\xi_1$ imply $b^\ast\leq h$, and then also $h^\ast\leq b$, which is, again, impossible because of $h^\ast\in\xi\cap\xi_1$.
    \end{itemize}
    \item \underline{Suppose $h^\ast\notin\xi_1$ and show $h\in\xi_1$:}

    Again, consider the two possibilities for $h$:
    \begin{itemize}
        \item If $h^\ast\in\xi$ then $h^\ast\leq b$, and consequently $b^\ast\leq h$ implying $h\in\xi_1$;
        \item If $h\in\xi$ but $h\notin\xi_1$, then $h\leq b$, which means $b^\ast\leq h^\ast$, implying $h^\ast\in\xi_1$ -- a contradiction.
    \end{itemize}
    \item \underline{Assume $h,k\in\xi_1$ with $h<k^\ast$, and derive a contradiction:}

    In the same manner as before, we conduct a case-study according to $h$ belonging or not belonging to $\xi$.
    \begin{itemize}
        \item If $h\in\xi$ then $k^\ast\in\xi$; since $k\in\xi_1$, we conclude $k^\ast\leq b$. In particular, $h\leq b$, so $h\notin\xi_1$ -- contradiction.
        \item Thus, $h^\ast\in\xi$ and $h\in\xi_1$, and we have $h^\ast\leq b$, and then also $k\leq b$. But $k\in\xi_1$ means $k\leq b$ is impossible, and we are done.\ep
    \end{itemize}
\end{itemize}

\subsection{The Roller boundary}\label{The Roller boundary defined} It seems that Roller (\cite{[Rol]}) has been the first to consider topological properties of $H^\circ$ in the context of a ``Stone-type'' duality he had discovered between median algebras and poc-sets. Therefore, in view of our cor.\ref{cor:ordering the Roller boundary}, we define 
\begin{defn}[Roller boundary of a poc-set] Let $(H,\leq,\ast)$ be an $\omega$-dimen\-sional dis\-crete poc-set. 

The \emph{Roller boundary} $\Re H$ of $H$ is the set of almost-equality classes of ultrafilters, partially-ordered by the relation
\begin{eqnarray}
    \Sigma_1\geq \Sigma_2 & \stackrel{\scriptscriptstyle{d\!e\!f\!\!.}}{\IFF} &
    \Sigma_1\cap\cl{\Sigma_2}\neq\varnothing\\
        &\IFF& \Sigma_1\subseteq\cl{\Sigma_2}.
\end{eqnarray}
\end{defn}

\subsubsection{Codimension}\label{subsubsection:codimension}
For our purposes it will be beneficial to develop a notion of (co)dimension in $\Re H$.

Let us fix $\xi\in H^\circ$ and denote its class with $\Sigma$. Assume $\Sigma$ is not the principal class.

Consider the set of all infinite descending chains in $H$. If $c=\{c_n\}_{n=1}^\infty$ and $d=\{d_m\}_{m=1}^\infty$, say that \emph{$c$ dominates $d$} (denoted $c\succ d$), if for every $m$ there exists $n$ such that $c_n\leq d_m$. We say the chains $c$ and $d$ are \emph{equivalent} (denoted by $c\sim d$), if both $c\succ d$ and $d\succ c$ hold. It is immediate that the relation $(\succ)$ induces a partial ordering on the set of equivalence classes of chains, and that every chain is equivalent to any of its infinite sub-chains.  

We apply the above relations to the set $D_\xi$ of infinite descending chains which are contained in $\xi$. Our main interest lies with the set $D_\xi/\!\sim$. 

Note that for any $\sigma\in\Sigma$ and any chain $c\in D_\xi$, $c$ is eventually-contained in $\sigma$, so that each class of $D_\xi$ corresponds to a unique class in $D_\sigma$. Thus, for example, the cardinality of $D_\xi/\!\sim$ is an invariant of $\Sigma$.
\begin{defn}[codimension of an a.e. class]\label{defn:codim} Let $\Sigma\in\Re H$ be any almost-equality class of $H^\circ$. The \emph{codimension of $\Sigma$ in $H$} is defined to be the cardinality $D_\xi/\!\sim$ for some (and hence any) $\xi\in\Sigma$.
\end{defn}
For example, the principal class is characterized by being of codimension zero.

Our current objective is to characterize the classes of finite codimension in terms of the order structure we have for $\Re H$.

Suppose now that $\xi\in\Sigma$, where $\Sigma\in\Re H$ is of finite non-zero codimension $\delta$. Now, $D_\xi/\!\sim$ carries the domination partial ordering, and we may select an equivalence class by fixing a corresponding chain $c=\{c_n\}_{n=1}^\infty$. For each $n$, we let $\xi_n$ be the result of truncating $\xi$ at $c_n$. 

We note the sequence $(\xi_n)_{n=1}^\infty$ converges in $H^\circ$ to the limit $\xi$, tempting us to wish for eventually-all of the $\xi_n$ to belong to the same almost-equality class. Were this our case, then we would have found a $\Sigma_1\in\Re H$ such that $\Sigma\cap\cl{\Sigma_1}$ is non-empty; by corollary \ref{cor:ordering the Roller boundary}, we would have $\Sigma_1<\Sigma$ in $\Re H$.

Let us prove that this is indeed the case when one selects the equivalence class of $c$ carefully enough: from now on, assume this class is maximal with respect to the domination ordering. The existence of such a class is guaranteed by our assumption that the codimension of $\Sigma$ is finite.

Fix a complete set of representatives $c^{(1)}=c,\ldots,c^{(\delta)}$ in $D_\xi/\sim$. For each $i>1$, we have $c^{(i)}\not\succeq c$, meaning that there exists $N_i\in\NN$ such that no relation of the form $c^{(i)}_m\leq c_{N_i}$ is possible for any $m$. We choose $N$ to be maximal among the $N_i$. 

Given $n>N$, $\xi_n$ may be obtained from $\xi_{n+1}$ by truncating $c_{n+1}$ (see lemma \ref{lemma:transitivity of truncation}); thus, in order to show these two ultrafilters are almost-equal, it will be enough to prove that $c_n$ does not lie on an infinite descending chain through $\xi_{n+1}$. 

Suppose $d=\{d_m\}_{m=1}^\infty$ is a descending chain in $\xi_{n+1}$ with the property that $d_1=c_n$. Since $d$ is a descending chain and $H$ is discrete, its elements have to eventually leave the set $\{h\in H\,\big|\,c^\ast_{n+1}\leq h\}$, so that a terminal segment of $d$ represents a chain-class in $\xi$. Now, $d$ cannot be equivalent to $c$, and we conclude that $c_n=d_1\geq c^{(i)}_m$ for some fixed $i$ and any $m$ large enough, contradicting the definition of $N$. 

One last remark has to be made regarding this construction: by the construction of $\xi_n$ (for $n>N$), it is evident that all the chains $c^{(i)}$ ($i>2$) survive in $\xi_n$, which means that by passing from $\Sigma$ to $\Sigma_1$ we have reduced the codimension of $\Sigma$ by precisely $1$. 

We summarize our last results in the following
\begin{prop}\label{prop:characterizing dimension} Suppose $\Sigma\in\Re H$ has finite dimension
$\delta>0$. Then, there exists $\Sigma_1\in\Re H$ such that:
\begin{itemize}
    \item[$(i)$] $\Sigma_1<\Sigma$, and
    \item[$(ii)$] $codim(\Sigma_1)=\delta-1$.
\end{itemize}
\end{prop}

In particular, we have shown that $\Sigma$ has codimension $\delta$ if and only if the shortest maximal descending chain in $\Re H$ starting at $\Sigma$ (and terminating in $\Pi$) has length $\delta$.

\subsubsection{Natural median structure on $\Re H$.} We recall $H^\circ$ has the structure of a median algebra defined by the median operation 
\begin{equation}
    \med(\alpha,\beta,\gamma)=(\alpha\cap\beta)\cup(\alpha\cap\gamma)\cup(\beta\cap\gamma).
\end{equation}
This operation is well-defined modulo almost-equality, inducing a median operation (satisfying all the required identities) on $\Re H$. It will be important to us that this operation is, in some sense, well-behaved with respect to the ordering we have defined on $\Re H$:
\begin{lemma}[principal intervals]\label{lemma:principal intervals} Suppose $\Sigma\in\Re H$, and $\Pi$ is the principal class in $H^\circ$. Then
\begin{equation}
    [\Pi,\Sigma]_{med}=\left\{\Sigma_1\,\big|\,\Sigma_1\leq\Sigma\right\},
\end{equation}
where $[\cdot,\cdot]_{med}$ refers to the interval in $\Re H$ with respect to the median operation.
\end{lemma}
\proof{} It follows from the properties of median functions that the interval $[\Pi,\Sigma]_{med}$ may be defined as (compare with the definition in the beginning of \ref{convexity structure on H circ}) 
\begin{equation}
    [\Pi,\Sigma]_{med}=\left\{\med(\Pi,\Sigma,A)\,\big|\,A\in\Re H\right\}.
\end{equation}
Suppose now $A\in\Re H$, and consider $\Sigma_1=\med(\Pi,\Sigma,A)$. We need to show two things: that $\Sigma_1\leq\Sigma$, and that every $B\in\Re H$ with $B\leq\Sigma$ arises in the same manner as did $\Sigma_1$. 

To show that $\Sigma_1\leq\Sigma$ it is enough to produce an element of $\Sigma$ lying in the Tychonoff closure of $\Sigma_1$. Fixing $\alpha\in A$, $\sigma\in\Sigma$ we consider the ultrafilters 
\begin{equation}
\sigma_\pi=\med(\pi,\sigma,\alpha)=(\alpha\cap\sigma)\cup(\alpha\cap\pi)\cup(\sigma\cap\pi)
\end{equation}
in $\Sigma_1$ arising for different choices of $\pi\in\Pi$. Given a finite subset $F$ of $\sigma$, use the denseness of $\Pi$ in $H^\circ$ to find $\pi\in\Pi$ containing $F$. Then the formula for $\sigma_\pi$ shows $F\subset\sigma_\pi$, proving $\sigma$ lies in the Tychonoff closure of $\Sigma_1$, as desired.

Conversely, suppose $\Sigma_1\leq\Sigma$, let $\Sigma_2=\med(\Pi,\Sigma,\Sigma_1)$, and let us show $\Sigma_2=\Sigma_1$. By what we have shown already, using the symmetries of the median operation, we see that both $\Sigma_2\leq\Sigma$ and $\Sigma_2\leq\Sigma_1$.

Therefore, it will be enough to prove $\Sigma_1\leq\Sigma_2$ -- that is, to show $\Sigma_2$ is contained in the Tychonoff closure of $\Sigma_1$. Consider an element $\sigma_2\in\Sigma_2$ arising in the form
\begin{equation}
    \sigma_2=(\sigma_1\cap\sigma)\cup(\sigma_1\cap\pi)\cup(\sigma\cap\pi),
\end{equation}
where $\pi\in\Pi$, $\sigma\in\Sigma$ and $\sigma_1\in\Sigma_1$. If now $F$ is a finite subset of $\sigma_2$ and $\pi\in\Pi$, as before, is a principal ultrafilter containing $F$, then our aim is to show that
\begin{equation}
    \Sigma_1\cap\bigcap_{f\in F}S_f\neq\varnothing.
\end{equation}
Since $F$, being a subset of $\sigma_2$, is a filter base in $H$, our version of Helly's theorem (see corollary\ref{cor:Helly's theorem}) implies it is enough to show $\Sigma_1\cap S_f$ is non-empty for every $f\in F$.

We finish the argument by way of contradiction. Suppose some $f\in F$ has $\Sigma_1\subset S_f^\ast$. Since $S_f^\ast$ is Tychonoff-closed, $\cl{\Sigma_1}$ -- and therefore $\Sigma$, too -- is contained in $S_f^\ast$. In particular, neither of the sets $\sigma_1\cap\sigma$, $\sigma_1\cap\pi$ and $\sigma\cap\pi$ contains $f$, implying $\sigma_2$ does not contain it -- a
contradiction.\ep\\

We shall henceforth supress the subscript $_{med}$ in the notation for intervals in $\Re H$, knowing that both notions of an interval in $\Re H$ -- both as a median algebra and as a partially ordered set -- coincide. 
\begin{defn}[A notion of gcd in $\Re H$]\label{defn:gcd} For $A,B,\Sigma\in\Re H$, we say $\Sigma$ is a greatest common divisor for $A$ and $B$, and denote this by $\Sigma=gcd(A,B)$, if
\begin{enumerate}
    \item $\Sigma\leq A$ and $\Sigma\leq B$, and --
    \item for any $\Sigma'\in\Re H$ satisfying $\Sigma'\leq A$ and $\Sigma'\leq B$ we have $\Sigma'\leq\Sigma$.
\end{enumerate}
\end{defn}
Greatest common divisors indeed exist, as the following corollary of the last lemma shows:
\begin{cor}[existence of a gcd]\label{cor:gcd's exist} Any two elements $A,B\in\Re H$ have a gcd. Moreover, 
\begin{equation}
    gcd(A,B)=\med(\Pi,A,B),
\end{equation}
where $\Pi$ is the principal class of $H$.
\end{cor}
\proof{} Set $\Sigma=\med(A,B,\Pi)$ and check it satisfies the required conditions.
\begin{enumerate}
    \item by definition, $\Sigma\in[\Pi,A]$, and the last lemma implies $\Sigma\leq A$. The same argument works for $B$.
    \item Since $\Sigma'\leq A,B$, we have, by the last lemma,
    \begin{eqnarray}
        \Sigma'&\in&[\Pi,A]\cap[\Pi,B]\\
        &=&[\Pi,\med(\Pi,A,B)]
    \end{eqnarray}
    (for the equality, see \cite{[Rol]}, pg.15 claim \textbf{(Int 5)}), and use the same lemma again to deduce $\Sigma'\leq\Sigma$.\ep
\end{enumerate}

\subsubsection{Restriction to a sub-poc-set.}\label{subsubsection:restriction on boundaries} Suppose $K$ is a sub-poc-set of the $\omega$-dimensional poc-set $H$ (i.e., $K$ is closed under $h\mapsto h^\ast$ and contains $0$). We then have the restriction map $res$ arising as the dual of the inclusion of $K$ into $H$ and preserving almost-equality. Consequently, we obtain a well-defined map $R^H_K:\Re H\to\Re K$. 
\begin{lemma}\label{lemma:restriction is monotone} $R^H_K$ is a monotone non-decreasing map, i.e., if $\Sigma_1\leq\Sigma_2$ in $\Re H$, then $R^H_K\Sigma_1\leq R^H_K\Sigma_2$ in $K$.
\end{lemma}
\proof{} Saying $\Sigma_1\leq\Sigma_2$ is equivalent to saying $\Sigma_2\subset\cl{\Sigma_1}$. Since $res$ is a continuous map, we conclude $res(\Sigma_2)\subseteq \cl{res(\Sigma_1)}$. Thus, the Tychonoff closure of $R^H_K\Sigma_1$ intersects $R^H_K(\Sigma_2)$, showing the desired inequality.\ep

\subsection{Convergence of principal ultrafliters.}
Let us now consider some issues of convergence (in the Tychonoff topology) of sequences of (principal) ultrafilters.
\begin{lemma}\label{lemma:joint convergence of UF's} Suppose $(\sigma_n)_{n=1}^\infty$, $(\tau_n)_{n=1}^\infty$ are convergent sequences in $H^\circ$ with limits $\sigma$ and $\tau$, respectively. If there exists $M\in\NN$ such that $\Delta(\sigma_n,\tau_n)\leq M$ for all $n$, then $\sigma\ae\tau$.
\end{lemma}
\proof{} Suppose $\sigma\minus\tau$ contains an infinite sequence $\{h_m\}_{m=1}^\infty$, and let us prove $\Delta(\sigma_n,\tau_n)$ is unbounded. Note that $h_m^\ast\in\tau$ for all $m$.

Given $m\in\NN$ we define
\begin{equation}
    B_m=\left\{h_1,\ldots,h_m\right\},
\end{equation}
producing a pair of disjoint open neighbourhoods $V(B_m)$ and $V(B_m^\ast)$ of $\sigma$ and $\tau$, respectively. By the definition of convergence in $H^\circ$, there exists $N$ such that $\sigma_n\in V(B_m)$ and $\tau_n\in V(B_m)$ hold for all $n\geq N$, showing $\Delta(\sigma_n,\tau_n)$ is unbounded.\ep\\

The next result, shows how one can ``optimize'' convergent sequences of principal ultrafliters.
\begin{lemma}[averaging lemma]\label{lemma:averaging lemma} Suppose $(\pi_n)_{n=1}^\infty$ is an infinite sequence of principal ultrafilters converging to the limit $\sigma\in H^\circ$. Let $\Sigma$ denote the almost-equality class of $\sigma$, and let $pr_\Sigma$ denote the projection map of $\Pi$ onto $\Sigma$. 

We define a new sequence of principal ultrafilters inductively as follows:
\begin{eqnarray}
    \pi'_1&=&\pi_1\;,\\
    \pi'_2&=&\med(\pi'_1,\pi_2,\sigma)\;,\\
    &\vdots&\\
    \pi'_{n+1}&=&\med(\pi'_n,\pi_{n+1},\sigma)\;.
\end{eqnarray}
Then, the sequence we have defined converges on $\sigma$ and satisfies the additional requirement that $\Delta(pr_\Sigma(\pi'_n),\sigma)$ is an eventually-zero monotone non-increasing function of $n$. 
\end{lemma}
\proof{} Suppose now $A$ is a finite subset of $\sigma$. Find $N$ such that
$\pi_n\in V(A)$ for all $n\geq N$. But then, for all $n\geq N$ we
must have $A\subseteq\pi_n\cap\sigma\subseteq\pi'_n$, showing
$\pi'_n\in V(A)$, and proving that the $\pi'_n$ do converge on
$\sigma$.

Now let us compare $\pi'_{n+1}\vartriangle\sigma$ to $\pi'_n\vartriangle\sigma$:
\begin{eqnarray}
    h\in\sigma\minus\pi'_{n+1}&\IFF& h\in\sigma\;\wedge\;h^\ast\in\pi'_{n+1}\\
    &\THEN& h\in\sigma\;\wedge\;h^\ast\in\pi'_n\cap\pi_{n+1}\\
    &\THEN& h\in\sigma\;\wedge\;h^\ast\in\pi'_n\\
    &\IFF& h\in\sigma\minus\pi'_n\;.
\end{eqnarray}
Thus, $\sigma\minus\pi'_{n+1}\subseteq\sigma\minus\pi'_n$. Applying the complementation operator to both sides we obtain $\pi'_{n+1}\minus\sigma\subseteq\pi'_n\minus\sigma$, producing the inequality
\begin{equation}
	\pi'_{n+1}\vartriangle\sigma\subseteq\pi'_n\vartriangle\sigma.
\end{equation}
The result of applying $pr_\Sigma$ to both sides of the inequality is obtained by intersecting both sides of the inequality with $T(\Sigma)$ (and then applying the natural isomorphism between $\cl{\Sigma}$ and $T(\Sigma)^\circ$), and we obtain:
\begin{equation}
	pr_\Sigma(\pi'_{n+1})\vartriangle\sigma\subseteq
	pr_\Sigma(\pi'_n)\vartriangle\sigma
\end{equation}
This shows that the quantity $\Delta(pr_\Sigma(\pi'_n),\sigma)$ is, indeed, a non-increasing function of $n$ taking values in $\NN\cup\{0\}$. Moreover, the above inclusion shows that, if $v$ is the eventual value of the latter quantity, then the above inequality implies there is an element $\sigma_1\in\Sigma$ which equals eventually all the projections $pr_\Sigma(\pi'_n)$. Therefore, $\sigma_1$ is a limit of the sequence $pr_\Sigma(\pi'_n)$. Since $pr_\Sigma$ is continuous and the $\pi'_n$ converge to $\sigma$, the fact $H^\circ$ is Hausdorff implies $\sigma=\sigma_1$, showing $v=0$, as desired.\ep\\

We now consider a special case of convergence.
\begin{lemma}[geodesic sequences converge]\label{lemma:geodesic sequences converge} Suppose $(\pi_n)_{n=1}^\infty$ is a geodesic ray in $\Pi$ with respect to the metric $d$ -- that is, for all $m,n\in\NN$ we have
\begin{equation}
    \Delta(\pi_n,\pi_m)=\big|n-m\big|.
\end{equation}
Then $(\pi_n)_{n=1}^\infty$ converges to a limit $\sigma\in H^\circ$, and the averaged sequence $(\pi'_n)_{n=1}^\infty$ defined above coincides with $(\pi_n)_{n=1}^\infty$. 
\end{lemma}
\proof{} We may write $\pi_{n+1}=[\pi_n]_{a_n}$ for $a_n\in\min(\pi_n)$ and every $n\in\NN$.  The assumption that the $\pi_n$ form a geodesic ray in $\Pi$ means that 
\begin{equation}
    \pi_i\vartriangle\pi_j=\left\{a_j,\ldots,a_{i-1},a^\ast_j,\ldots,a^\ast_{i-1}\right\}
\end{equation}
for all $i>j$, with all the $a_n$ distinct. In particular, for $j=1,i=n>1$ we have that the list
\begin{equation}
    L_n:=\pi_{n+1}\minus\pi_1=\left\{a^\ast_1,\ldots,a^\ast_n\right\},
\end{equation}
consists of distinct $n$ distinct elements and that $\pi_m$ contains $L_n$ for all $n>m$.

Let now $h\in H$ be a fixed halfspace. We shall show that either $h\in\pi_n$ for all but finitely many values of $n$, or $h^\ast\in\pi_n$ for eventually all $n$: indeed, if $h\notin\pi_n$ for all $n$, then $h^\ast\in\pi_n$ for all $n$ and we are done; the same holds for $h^\ast$, so we may assume there are $i,j\in\NN$ such that $h\in\pi_i$ and $h^\ast\in\pi_j$; without loss of generality, assume $i<j$, and we are forced to conclude that $h^\ast$ appears on the list $L_j$, meaning $h$ cannot appear on any of the lists $L_n$ for $n>j$. Thus, $h^\ast\in\pi_n$ for all $n>j$, as desired.

We define an ultrafilter $\sigma$ as follows: for any $h\in H$ we set $h\in\sigma$ if and only if $h\in\pi_n$ for all but finitely-many values of $n$. By what we have just shown, $\sigma$ is, indeed, an ultrafilter on $H$, and since $(\pi_n)_{n=1}^\infty$ converges to $\sigma$ pointwise, it converges to $\sigma$ in the Tychonoff topology, as required.\\

Finally, let us use induction to prove $\pi'_n=\pi_n$ holds for all $n$, the base step ($n=1$) having been taken care of by the definition of the ``averaged'' sequence. Given $n\in\NN$, suppose that for all $k\leq n$ we have $\pi'_k=\pi_k$, and let us consider $\pi'_{n+1}$:
\begin{eqnarray}
    \pi'_{n+1}&=&(\pi'_n\cap\pi_{n+1})\cup(\pi'_n\cap\sigma)\cup(\pi_{n+1}\cap\sigma)\\
    &=&(\pi_n\cap\pi_{n+1})\cup(\pi_n\cap\sigma)\cup(\pi_{n+1}\cap\sigma)\\
    &=&(\pi_{n+1}\minus\{a^\ast_n\})\cup(\pi_n\cap\sigma)\cup(\pi_{n+1}\cap\sigma).
\end{eqnarray}
Note that the second and third summands differ only by the element $a^\ast_n$, which, by the argument given above (set $h=a_n$) must lie in $\sigma$. Thus, $\pi'_{n+1}$ contains $\pi_{n+1}$ as a subset, implying they are equal, as desired.\ep
\begin{cor}\label{cor:geodesic sequences converge tamely} Suppose $(\pi_n)_{n=1}^\infty$ is a geodesic ray in $\Pi$, and let $\sigma$ be its limit. Further, let $\Sigma$ be the almost-equality class of $\sigma$, and let $pr_\Sigma$ denote the natural projection map of $\Pi$ onto $\Sigma$. Then $pr_\Sigma(\pi_n)$ is eventually constant, and $\Delta(pr_\Sigma(\pi_n),\sigma)$ is monotone non-increasing and eventually-zero.\ep
\end{cor}

The previous result has a converse. It is only natural to ask how far is an averaged sequence from being a geodesic vertex-path in $\Pi$:
\begin{prop}\label{prop:averaged is subsequence of geodesic} Suppose $(\pi_n)_{n=1}^\infty$ is a sequence in $\Pi$ converging to a limit $\pi_\infty$ in $H^\circ$, and let $(\pi'_n)_{n=1}^\infty$ denote the averaged sequence as defined above. Then $(\pi'_n)_{n=1}^\infty$ is a subsequence of a geodesic ray in $\Pi$.
\end{prop}
\proof{} Given the elements $\pi'_n$, we select, for each $n\geq 1$, a geodesic vertex path $p_n$ from $\pi'_n$ to $\pi'_{n+1}$. Our guess is that the concatenation of all these paths is a geodesic ray in $\Pi$. 

Let $q_n$ denote the concatenation $p_1\ast\ldots\ast p_n$. It will be enough to show that $q_n$ is a geodesic path for all $n$. To do that, we use induction on $n$, where the case $n=1$ is a sure success, by construction.

Suppose now that $q_{n-1}$ is a geodesic vertex path in $\Pi$ for some $n\geq 2$. In view of the induction hypothesis it is enough to prove the following equality: 
\begin{displaymath}
	\Delta(\pi'_1,\pi'_n)+\Delta(\pi'_n,\pi'_{n+1})=\Delta(\pi'_1,\pi'_{n+1}).
\end{displaymath}
Since $q_{n-1}$ is a geodesic, to prove this equality it would be enough to show that
\begin{displaymath}
	\Delta(\pi'_{n-1},\pi'_n)+\Delta(\pi'_n,\pi'_{n+1})=\Delta(\pi'_{n-1},\pi'_{n+1}).
\end{displaymath}
Applying the definition of $d$, this takes the form of:
\begin{displaymath}
	\big|\pi'_{n-1}\vartriangle\pi'_n\big|+\big|\pi'_n\vartriangle\pi'_{n+1}\big|=
	\big|\pi'_{n-1}\vartriangle\pi'_{n+1}\big|.
\end{displaymath}
Since $\pi'_{n-1}\vartriangle\pi'_{n+1}=(\pi'_{n-1}\vartriangle\pi'_n)\vartriangle
(\pi'_n\vartriangle\pi'_{n+1})$, proving this would be the same as to prove that $\pi'_{n-1}\vartriangle\pi'_n$ and $\pi'_n\vartriangle\pi'_{n+1}$ do not intersect.

Note both sets are symmetric, and hence so is their intersection. 
Suppose now there was an element $h\in\pi'_{n+1}\minus\pi'_n$ lying in $\pi'_n\vartriangle\pi'_{n-1}$. Observe that then $h\in\pi'_{n-1}\minus\pi'_n$.
Since
\begin{displaymath}
	\pi'_{n+1}=(\pi'_n\cap\pi_{n+1})\cup(\pi'_n\cap\pi_\infty)\cup
	(\pi_{n+1}\cap\pi_\infty),
\end{displaymath}
and since $h\notin\pi'_n$, we conclude $h\in\pi_{n+1}\cap\pi_\infty$.

On the other hand, writing
\begin{displaymath}
	\pi'_n=(\pi'_{n-1}\cap\pi_n)\cup(\pi'_{n-1}\cap\pi_\infty)\cup
	(\pi_n\cap\pi_\infty),
\end{displaymath}
and observing $h^\ast\in\pi'_n\minus\pi'_{n-1}$ imply $h^\ast\in\pi_n\cap\pi_\infty$.

Thus, both $h$ and $h^\ast$ lie in $\pi_\infty$ -- a contradiction to $\pi_\infty$ being an ultrafilter.\ep\\

In order to complete the picture it remains to show a standard way of approximating a given element $\sigma_\infty$ of a given class $\Sigma\in\Re H$ by an averaged sequence of elements of $\Pi$.

For every $\pi\in\Pi$, recall that $pr_\Sigma(\pi)=P(\Sigma)\cup\left(\pi\cap T(\Sigma)\right)$ is an element of $\Sigma$. In fact, we know that $pr_\Sigma:\Pi\to\Sigma$ is a surjective median morphism, which allows to select $\pi\in pr_\Sigma^{-1}(\sigma_\infty)$. 

Consider the intersection of $\min(\pi)$ with $P(\Sigma)^\ast$: for any $a,b\in\min(\pi)$, the only possible relations are $a\pitchfork b$ and $a^\ast\leq b$; however, if $\Sigma\subset S_{a^\ast}$, then $a^\ast\leq b$ would imply $\Sigma\subset S_b$, not allowing $\Sigma\subset S_{b^\ast}$. Thus, $\min(\pi)\cap P(\Sigma)^\ast$ is a transverse subset of $\min(\pi)$. We may then define a function 
\begin{equation}
	\left\{\begin{array}{rcl}
		F_\Sigma:\Pi&\to&\Pi\\
		\pi&\mapsto&\left[\pi\right]_{\min(\pi)\cap P(\Sigma)^\ast}
	\end{array}\right.
\end{equation}
The sequence $\sigma_n=F_\Sigma^n(\pi)$ then necessarily converges to $pr_\Sigma(\pi)$ in the Tychonoff topology. Moreover, $(\sigma_n)_{n=1}^\infty$ is easily seen to be an averaged sequence once we notice that $\sigma_n\cap T(\Sigma)$ is constant and that the intersections $\sigma_n\cap P(\Sigma)$ form a strictly increasing sequence.\\

This construction is motivated by the notion of a normal cube path introduced by Niblo and Reeves in \cite{[NibRee2]}, with the $\sigma_n$ constituting, in a sense, a normal cube path from $\pi$ to $pr_\Sigma(\pi)$.
\begin{defn}[canonical flow]\label{defn:canonical flow} Given $\Sigma\in\Re\HH$, the map $F_\Sigma:\Pi\to\Pi$ will be called {\it the canonical flow on $\Pi$ in the direction of $\Sigma$}. 
\end{defn}
The motivation for the name is, of course, the fact that iterating $F_\Sigma$ over $\Pi$ decomposes $\Pi$ into orbits converging (through iterating $F_\Sigma$) onto the points of $\Sigma$.

\section{Halfspace systems and boundary decompositions.}\label{section:boundary decompositions} Suppose $X$ is a proper $CAT(0)$-space.
\begin{defn} A \emph{halfspace} $h\subset X$ is a non-empty open convex subset such that $h^\ast\triangleq X\minus\cl{h}$ is also convex. The intersection $\cl{h}\cap\cl{h^\ast}$ will be called \emph{the wall associated with $h$}, and denoted by $W(h)$; if $S$ is a set of halfspaces, then $W(S)$ will denote the set of walls $W(h)$ for $h\in S$. The sets $\varnothing, X$ are, by definition, the \emph{trivial} halfspaces of $X$. 
\end{defn}
Our notion of a halfspace system is the obvious generalization of what one observes in a cubing. 
\begin{defn} A \emph{halfspace system} in a CAT(0) space $X$ is a family $\HH$ of
halfspaces containing the trivial halfspaces, ordered by containment, invariant under the operation $h\mapsto h^\ast$ and satisfying the following conditions: 
\begin{itemize}
    \item[$\mathbf{(H1)}$] Every point $x\in X$ has a neighbourhood $U_x$ intersecting only finitely many walls associated with halfspaces of $\HH$.
    \item[$\mathbf{(H2)}$] $\HH$ contains no infinite transverse subfamily. 
\end{itemize}
\end{defn}
The above notion of a halfspace system is discussed in detail in \cite{[Gu-loc-fin]}. The additional requirement $(H2)$ is required for us to be able to apply the results of the previous section. For the sake of completeness, we provide all the definitions in this text, but we shall refer the reader there for the proofs of basic results.

The goal of this section is to explain how the Roller boundary $\Re\HH$ of a halfspace system $\HH$ induces a partition on the CAT(0)-boundary $\bd X$ of $X$, and to study some of the simplest connections between geometric properties of $\HH$ and the behavior of $\rho$ in order to establish a reasonably general framework within which to continue our discussion.

Two facts about geodesic rays and halfspaces in CAT(0) spaces are used throughout this work without reference:
\begin{lemma}\label{lemma:rays and halfspaces1} Suppose $h$ is a halfspace in $X$, and suppose $c:[0,\infty)\to X$ is a geodesic ray  emanating from a point of $h$. If $c(\infty)\in\bd h$, then the image of $c$ lies in the closure of $h$. Moreover, if $c(T)\in W(h)$ for some $T>0$, then $c(t)\in W(h)$ for all $t\geq T$, and $c(\infty)\in\bd W(h)$.\ep
\end{lemma}
\begin{lemma}\label{lemma:rays and halfspaces2} For any halfspace $h$ in $X$ we have $\bd W(h)=\bd h\cap\bd h^\ast$. \ep
\end{lemma}
%
%
%
%

\subsection{Defining the boundary decomposition map.}

Given a halfspace system $\HH$ on $X$ and a point $\xi\in\bd X$ we consider the set 
\begin{equation}
	\HH^\xi=\left\{
		\sigma\in\HH^\circ\,\bigg|\,
		\forall_{h\in\sigma}\;\xi\in\bd\cl{h}
	\right\}.
\end{equation}
Given any $\sigma\in\HH^\xi$, the set 
\begin{equation} 
	T(\xi)=\left\{
		h\in\HH\,\bigg|\,
		\xi\in\bd\cl{h}\minus\bd W(h)
	\right\}
\end{equation}
is obviously contained in $\sigma$. One should think of $W(h)$ as being transverse to $\xi$ in the sense that for any point $p\in h^\ast$ the geodesic ray $\gamma$ emanating from $p$ and having endpoint $\xi$ will inevitably cross $W(h)$ -- hence the notation for $T(\xi)$. An immediate observation following from the preceding lemmas is:
\begin{lemma}\label{lemma:when a bdry point lies in a filter} Let $\sigma\in\HH^\circ$ and $\xi\in\bd X$. Then $\sigma\in\HH^\xi$ if and only if $\sigma$ contains $T(\xi)$.\ep
\end{lemma}
\proof{} We had already noticed that $T(\xi)\subseteq\sigma$ for all $\sigma\in\HH^\xi$. Suppose now that $\sigma\in\HH^\circ$ contains $T(\xi)$, and let us show it lies in $\HH^\xi$. Thus, taking $h\in\sigma\minus T(\xi)$ we must show $\xi\in\bd\cl{h}$. Since $\sigma$ contains $T(\xi)$, we must have $h^\ast\notin T(\xi)$, and we have that neither of 
\begin{equation}
    \xi\in\bd\cl{h}\minus\bd W(h)\,,\qquad\xi\in\bd\cl{h^\ast}\minus\bd W(h)
\end{equation}
holds. This implies $\xi\in\bd W(h)$, finishing the proof.\ep 
\begin{cor} $\HH^\xi$ is non-empty for all $\xi\in\bd X$.
\end{cor}
\proof{} $T(\xi)$ is clearly a filter base, and is therefore contained in an ultrafilter.\ep\\
 
Continuing the analysis of $\HH^\xi$ let us define
\begin{equation} 
    P(\xi)=\{0,0^\ast\}\cup\left\{h\in\HH\,\bigg|\,\xi\in\bd W(h)\right\}.
\end{equation}
Observe that $P(\xi)$ inherits the structure of a discrete $\omega$-dimensional poc-set. By defintion, $\sigma\in\HH^\xi$ implies $\sigma\subset T(\xi)\cup P(\xi)$. Thus, a halfspace $a\in\min(\sigma)$ will lie either in $T(\xi)$ or in $P(\xi)$:
\begin{itemize}
	\item[-] if $h\in\min(\sigma)\cap P(\xi)$, then $[\sigma]_h$ remains in $\HH^\xi$;
	\item[-] if $h\in\min(\sigma)\cap T(\xi)$, we inevitably have $[\sigma]_h\notin\HH^\xi$.
\end{itemize}
The first observation hints at a possible connection between the structure of $\HH^\xi$ and that of $P(\xi)^\circ$, realized through the map $res_\xi$ dual to the inclusion $i:P(\xi)\hookrightarrow\HH$. We shall now study the restriction of $res_\xi$ to $\HH^\xi$ more closely. 

For any $h\in P(\xi)$ and $k\in T(\xi)$, we have that $\xi$ is an interior point of $k$ and a boundary point of $h$ at the same time, which means $k<h$ is impossible; the relation $k<h^\ast$ is impossible for the same reason, as $h^\ast\in P(\xi)$. Thus, for all $h\in P(\xi)$ and $\sigma\in\HH^\xi$ we have 
\begin{equation}
	h\in\min(\sigma) \IFF h\in\min\left(\sigma\cap P(\xi)\right)=\min(res_\xi\sigma).
\end{equation}
We conclude that for every almost-equality class $\Sigma\in\Re\HH$ intersecting $\HH^\xi$, the map $res_\xi$ maps  $\Sigma\cap\HH^\xi$ \emph{onto} the corresponding almost-equality class of $P(\xi)^\circ$. 

Next, suppose $\bar\sigma\in P(\xi)^\circ$, and let us find a $\sigma\in\HH^\xi$ satisfying $\bar\sigma=res_\xi(\sigma)$. We shall use our previous observation that no element of $T(\xi)$ is smaller than an element of $P(\xi)$. Define $\sigma=T(\xi)\cup\bar\sigma$, and observe it automatically satisfies $\mathbf{(UF1)}$. In order to check $\mathbf{(UF2)}$ it is enough to take $h\in\bar\sigma$ and $k\in T(\xi)$ and rule out $h<k^\ast$ and $k<h^\ast$: applying the involution we see that the two are equivalent, and the second contradicts the observation we have just made. Thus we have proved:
\begin{prop}\label{prop:boundary decomposition map is well defined} Suppose $\HH$ is a halfspace system in a proper CAT(0)-space $X$, let $\xi\in\bd X$ be any point, and let $\Re(\xi)$ denote the set of almost-equality classes represented in $\HH^\xi$. Then, $res_\xi$ induces a one-to-one order-preserving correspondence $R_\xi$ of $\Re(\xi)$ onto $\Re P(\xi)$, with $res_\xi$ mapping $\Sigma\cap\HH^\xi$ bijectively onto $R_\xi(\Sigma)$. In particular, $\Re(\xi)$ has a minimum element -- the unique $\Sigma\in\Re(\xi)$ satisfying $R_\xi(\Sigma)=\Pi_{P(\xi)}$.\ep
\end{prop}
Thus the last proposition enables us to define a map of $\bd X$ into the Roller boundary of $\HH$:
\begin{defn}[boundary decomposition map] Suppose $\HH$ is a halfspace system in a proper CAT(0)-space $X$. The \emph{boundary decomposition map} 
\begin{equation}
	\rho:\bd X\to\Re\HH
\end{equation} 
is defined to be the map sending each $\xi\in\bd X$ to the minimum element of $\Re(\xi)$.
\end{defn}
The following examples provide the basic feeling of how this map works.
\begin{example} Let $\HH$ be the halfspace system arising from the standard cubing of Euclidean $2$-space and denote by $x^\pm$ and $y^\pm$ the positive and negative boundary points of the $x-$ and $y-$axes, respectively. For any $\xi\in\bd\EE^2$ other than $x^\pm,y^\pm$ and every  $h\in\HH$ we have either $h\in T(\xi)$ or $h^\ast\in T(\xi)$, showing $\HH^\xi$ consists of a unique point, and one easily sees there are precisely four such points; on the other hand, for $\xi=x^+$, say, one has $h\in T(\xi)$ for all positive vertical $h\in\HH$, whereas all horizontal elements of $\HH$ lie in $P(\xi)$. Thus, $H^\xi$ is homeomorphic to the ultrafilter space of the standard cubing of $\RR$, and $\rho$ maps the boundary point $x^+$ to the ``positive horizontal'' almost-equality class of $\HH^\circ$.
\end{example}
Note that if $\HH$ in the above example were taken to consist only of the horizontal halfspaces of the standard cubing of $\EE^2$, we would have all points $\xi$ of $\bd\EE^2$ except $y^\pm$ accept the principal class as their value under $\rho$, while $y^\pm$ would be mapped to the corresponding (two) non-principal classes. 
\begin{example} The same thing happens also in $X=\FAT{H}^2$, when $\HH$ is taken to be the family of halfspaces corresponding to the system of walls one obtains by translating a hyperbolic line $\ell$ by a hyperbolic translation $g$ with axis orthogonal to $\ell$. The line $\ell$ separates $g(+\infty)$ from $g(-\infty)$, making $\HH$ into a poc-set isomorphic to the standard poc-set structure on $\ZZ$, but here, given a point $\xi\in\bd X\minus\{g(\pm\infty)$ one is not able to associate a geodesic ray $\gamma$ satisfying $\gamma(\infty)=\xi$ with \emph{every} ultrafilter of the class $\rho(\xi)$.
\end{example}
Another interesting example is one of a `parabolic' nature:
\begin{example}
Again, let $X$ be the hyperbolic plane, and let $\ell$ be a line with an endpoint $p$ on the ideal boundary; denote by $h_\pm$ the different components of $X\minus\ell$. Further let $g$ be a parabolic translation fixing $p$, and let $\HH$ be the orbit of $\{h_\pm\}$ under $\langle g\rangle$. Then, for any ideal point $\xi\neq p$ we shall have $\rho(\xi)$ is the principal class $\Pi$ of $\HH^\circ$, with $\Re(\xi)$ containing only this class, while for $\xi=p$ the set $\Re(\xi)$ equals $\HH^\circ$, resulting in $\rho(\xi)=\Pi$ once again. Thus, here the representation map does not distinguish among the points of $\bd X$, though the assignment $\xi\mapsto\Re(\xi)$ does.
\end{example}
An immediate application of the construction of $\rho$ explains some of the behaviour of the above examples.
\begin{prop}\label{prop:when rho(xi)=Pi} A point $\xi\in\bd X$ satisfies $\rho(\xi)=\Pi$ if and only if $T(\xi)$ contains no infinite descending chains.
\end{prop}
\proof{} If $\rho(\xi)=\Pi$, then there exists $\sigma\in\Pi\cap\HH^\xi$. Since $T(\xi)\subseteq\sigma$, this implies $T(\xi)$ cannot contain an infinite descending chain. Conversely, if $\rho(\xi)>\Pi$, then an element $\sigma\in\rho(\xi)\cap\HH^\xi$ has to contain an infinite descending chain $c$. However, this chain cannot contain a subchain of elements from $P(\xi)$, as $\rho(\xi)\cap\HH^\xi$ is mapped by the restriction map into the principal class of $P(\xi)^\ast$. Since $\sigma\subset T(\xi)\cup P(\xi)$, all but finitely many elements of the chain $c$ must lie in $T(\xi)$.\ep

\subsection{Uniform systems.}
\begin{figure}[htb]
	\centering{\includegraphics[width=3in]{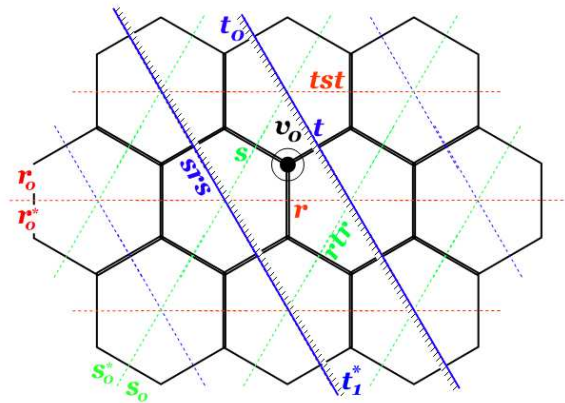}}
	\caption{\protect\scriptsize The hexagonal packing of $\EE^2$: walls are the fixed-point sets of reflections; $\HH$ decomposes as the union of three systems of proper halfspaces $\{r_n\},\{s_n\},\{t_n\}$ (and complements) indexed by $n\in\ZZ$, with $\{r_a,s_b,t_c\}$ transverse for any $a,b,c\in\ZZ$; we set $r_0,s_0,t_0$ to be three pairwise transverse minimal halfspaces among those containing the vertex $v_0$ corresponding to the unit element of $W$.\protect\normalsize}
	\label{figure:halfspace system}
\end{figure}

The motivation for our construction of $\rho$ was provided by known examples of cubings as well as by examples of tilings of the plane derived from the obvious Coxeter groups. We will henceforth restrict our attention to halfspace systems possessing a certain geometric property that is automatically satisfied when $X$ is a cubing and $\HH$ is its natural halfspace system. The same property is also enjoyed by the halfspace system defined on the Davis-Moussong complex of an infinite Coxeter system of finite rank: take $X$ to be the Davis-Moussong complex of a Coxeter system $(W,R)$ of finite rank, and let $\HH$ be the system of halfspaces arising as the set of complementary components of the walls. The example we have chosen for illustrating most of the work done in this paper is that of the regular hexagonal tiling of the Euclidean plane $\EE^2$, which is nothing else than the Davis-Moussong complex of the Coxeter system
\begin{equation}
	W\cong\left\langle r,s,t\left|r^2,s^2,t^2, (rs)^3,(rt)^3,(st)^3\right.\right\rangle\,,
\end{equation}
as illustrated in figure \ref{figure:halfspace system}.
Walls are defined to be the fixed point sets of reflections of the system $(W,R)$, and it can be shown (for example, see \cite{[Wil]}), that this system of walls coincides with $W(\HH)$.\\

The work of Brink and Howlett \cite{[BriHow]} shows that this particular class of examples has the \emph{parallel walls property}, which, in the case of Coxeter groups was introduced by Davis and Shapiro \cite{[Davis-Shapiro]}:
\begin{defn}[parallel walls property] A halfspace system $\HH$ in a proper CAT(0) space $X$ has the \emph{parallel walls property}, if there exists a constant $C>0$ such that for every $h\in\HH$ and $x\in X$ satisfying $d(x,h^\ast)>C$ there exists a halfspace $k\in\HH$ such that $x\in k<h$.
\end{defn}
\begin{remark} When we are mentioning parallel walls, this should not be mistaken for walls lying at a bounded Hausdorff distance from each other. By a pair of parallel walls we only mean walls $W(h),W(k)$ not being separated one from the other by another wall of $\HH$.
\end{remark}
With respect to the visual boundary $\bd X$ of $X$, $\HH$ having the parallel walls property becomes the following --
\begin{defn}[conical points, uniformness] A point $\xi\in\bd X$ is said to be a \emph{conical limit point of $\HH$}, $T(\xi)$ is non-empty, and for any $a\in T(\xi)$ and any cone neighbourhood $U$ of $\xi$ in $X\cup\bd X$ there exists $b\in T(\xi)$ satisfying $b<a$ and $b^\ast\cap U\neq\varnothing$.

A halfspace system $\HH$ on a proper CAT(0) space $X$ is said to be \emph{uniform}, if all points of $\bd X$ are conical limit points of $\HH$.
\end{defn}
\begin{example} A uniform halfspace system satisfies $\rho^{-1}(\Pi)=\varnothing$, by proposition \ref{prop:when rho(xi)=Pi}.
\end{example}
The following lemma in \cite{[Gu-loc-fin]} verifies that a halfspace system with the parallel walls property and satisfying $T(\xi)\neq\varnothing$ for all $\xi\in\bd X$ is uniform:
\begin{lemma}\label{lemma:T(xi) contains descending chains} Suppose $\HH$ is a halfspace system in a proper CAT(0) space $X$, having the parallel walls property, and let $\xi\in\bd X$. Then, for every $h\in T(\xi)$ and every cone neighbourhood $U$ of $\xi$ in $\hat{X}$ there exists a $k\in T(\xi)$ such that $k<h$ and $k^\ast\cap U\neq\varnothing$.\ep
\end{lemma}
Let us get back to the relation between our general situation and the situation in Coxeter groups. Given a halfspace system $\HH$, we notice that most points of the space do not lie on any wall of $\HH$. For any such point $x\in X$ it is possible to associate its {\it (closed) chamber} --
\begin{defn}[chamber] If $\HH$ is a halfspace system in a CAT(0) space $X$, and $x\in X$ does not lie on any wall of $\HH$, the chamber $ch(x)$ of $x$ is defined as the intersection of closures of all halfspaces in $\HH$ containing $x$.
\end{defn}
\begin{example} Clearly, since $\rho^{-1}(\Pi)=\varnothing$, a uniform halfspace system has bounded chambers.
\end{example}
In the Davis-Moussong $X=M(W,R)$ complex of a Coxeter system $(W,R)$, all chambers are bounded, as every chamber corresponds to a unique element of $W$, and $W$ acts co-compactly on $X$. The following result from \cite{[Gu-loc-fin]} establishes the connection between uniformness and the special case we had just considered:
\begin{thm} Let $X$ be a proper CAT(0) space and let $\HH$ be a halfspace system in $X$. If $\HH$ has bounded chambers and satisfies the parallel walls property, then $\HH$ is uniform.

If, in addition, $X$ admits a geometric group action by a group $G$ stabilizing $\HH$ and $\HH$ is uniform, then $\HH$ has the parallel walls property.\ep
\end{thm}

\subsubsection{Technical tools for uniform systems.}
We now derive some criteria for comparing the images of points $\xi,\eta\in\bd X$ under the decomposition map $\rho$ associated with a uniform system $\HH$.
\begin{lemma}\label{cor:comparison in bdd type} Suppose $\HH$ is a halfspace system in a proper CAT(0) space $X$, and let $\xi,\eta\in\bd X$. Then,
\begin{enumerate}
    \item If $T(\xi)\subseteq T(\eta)$, then $\rho(\xi)\leq\rho(\eta)$;
    \item If $\HH$ is uniform, then $\rho(\xi)\leq\rho(\eta)$ holds if and only if
    both $P(\eta)\subseteq P(\xi)$ and $T(\xi)\subseteq T(\eta)$.
\end{enumerate}
\end{lemma}
\proof{} Suppose $T(\xi)\subseteq T(\eta)$. Then $\rho(\eta)\in\Re(\xi)$, and since $\rho(\xi)$ is the minimum of $\rho(\xi)$, we get $\rho(\xi)\leq\rho(\eta)$.

Now, if $\HH$ is uniform, recall that for $\Sigma\in\Re\HH$, $T(\Sigma)$ denotes the poc-set of all $h\in\HH$ which are either trivial or satisfy $h\pitchfork\Sigma$ (i.e. $\Sigma\not\subset S_h$ and $\Sigma\not\subset S_h^\ast$). Consider $T(\rho(\xi))$: the uniformness of $\HH$ implies $T(\rho(\xi))=P(\xi)$ for all $\xi\in\bd X$.

Now assume $\rho(\xi)\leq\rho(\eta)$. Then $T(\rho(\xi))$ must contain $T(\rho(\eta))$, thus proving $P(\xi)$ contains $P(\eta)$. Next, recall $\rho(\xi)$ is contained in the closed set $\HH^\xi=\bigcap_{h\in T(\xi)}S_h$; since $\rho(\eta)$ is contained in the closure of $\rho(\xi)$, every element of $\rho(\eta)$ contains $T(\xi)$. Thus, $T(\eta)\cup T(\xi)$ is a filter base in $\HH$. Since $\HH=T(\eta)\cup T(\eta)^\ast\cup P(\eta)$ holds together with $T(\xi)\cap P(\eta)=\varnothing$, we conclude $T(\xi)\subset T(\eta)$.\ep
\begin{cor}\label{cor:how T(xi) grows moving upwards} Suppose $\rho(\xi)\leq\rho(\eta)$. Then
\begin{equation}
    T(\xi)\subseteq T(\eta)\subseteq T(\xi)\cup P(\xi).
\end{equation}
\end{cor}
\proof{} We have to show $h\in T(\eta)$ contradicts $h\in T(\xi)^\ast$. Take $\sigma\in\rho(\eta)$. Thus, $S_h$ is a neighbourhood of $\sigma$ and $\rho(\xi)\leq\rho(\eta)$ implies there exists $\tau\in\rho(\xi)$ lying in $S_h$, implying $h\in T(\xi)^\ast$ is impossible.\ep
\begin{lemma}\label{cor:comparison in bdd type - 2}
Suppose $\HH$ is a uniform halfspace system in a proper CAT(0) space $X$, and let $\xi,\eta\in\bd X$. If $T(\eta)\subseteq T(\xi)\cup P(\xi)$ and $P(\eta)\subseteq P(\xi)$, then $\rho(\xi)\leq\rho(\eta)$.
\end{lemma}
\proof{} Let us show $\HH^\eta\subset\HH^\xi$, which will imply $\rho(\xi)\leq\rho(\eta)$. Consider the restriction map $res_\eta:\HH^\eta\to P(\eta)^\circ$: we know it is an isomorphism of median algebras. If $\sigma\in\HH^\eta$, consider $\sigma\cap P(\xi)$: being a filter base, this set is contained in an element $\bar\sigma\in P(\xi)^\circ$. Set $\tilde\sigma=\bar\sigma\cup T(\xi)\in\HH^\xi$, and observe this is an ultrafilter containing all the elements of $\sigma$. Thus, $\sigma=\tilde\sigma$, and we have $\sigma\in\HH^\xi$.\ep

\subsubsection{The closure formula.}
\begin{prop}\label{prop:closure formula 0} Suppose $\HH$ is a uniform halfspace system in a proper CAT(0) space $X$, and let $\Sigma\in\Re\HH$ lie in the image of the boundary decomposition map $\rho$. Then --
\begin{equation}
\mathbf{(FF_0)}\quad\cl{\rho^{-1}(\Sigma)}\subseteq\bigcup_{\Sigma_1\leq\Sigma}\rho^{-1}(\Sigma_1).
\end{equation}
In particular, since $\rho^{-1}(\Pi)=\varnothing$, then $\rho^{-1}(\Sigma)$ is closed for all $\Sigma\in\Re\HH$ of codimension $1$.
\end{prop}
\proof{} Let $\xi\in\cl{\rho^{-1}(\Sigma)}$, and denote $\Sigma_1=\rho(\xi)$. In order to show $\Sigma_1\leq\Sigma$, it will be enough to show that $\Sigma$ is represented in $\HH^\xi$: since $\Sigma_1$ is a minimum in $\Re(\xi)$, we will have $\Sigma_1\leq\Sigma$. Thus, we must find $\sigma\in\Sigma$ containing $T(\xi)$.

Let $(\xi_n)_{n=1}^\infty$ be a sequence of boundary points satisfying $\rho(\xi_n)=\Sigma$ and converging to $\xi$. Recall that each of the $\HH^{\xi_n}$ is a union of almost-equality classes, and so one may fix an element $\sigma\in\Sigma$ and deduce that $T(\xi_n)\subseteq\sigma$ for all $n$. 

Now, if $h\in T(\xi)$, then there exists a neighbourhood $U$ of $\xi$ in $\bd X$  disjoint from $\bd\cl{h^\ast}$. Thus, $\xi_n\in U$ for $n$ large enough implies $h\in T(\xi_n)$, proving $h\in\sigma$.\ep\\

The above formula is, in a way, the main attribute of our decomposition map, demonstrating a relationship between the possibility of disconnecting large `chunks' (fibers of $\rho$) of $\bd X$ from each other and the structure of the image of $\bd X$ under $\rho$ as a partially-ordered set. We will see that more precision can be achieved for uniform systems.

\subsubsection{restriction to a subspace}\label{subsection:restriction to a subspace} Many examples of CAT(0) spaces arise as a result of various gluings. This is the reason why it might be a good startegy to compute the decomposition of the boundary of such a space based on the knowledge of the boundary decompositions of the subspaces being glued together to obtain the whole space. Proposition \ref{prop:the restriction equality} of this paragraph gives the precise tool for making such computations, showing that the boundary decomposition of any closed convex subspace of $X$ is the one induced from the decomposition of $X$ via the natural inclusion map.

Let us fix a closed convex subspace $F$ of $X$, and consider the map $r^X_F$ from $\HH$ into the set of halfspaces of $F$ defined by
\begin{equation}
    r^X_F(h)=\left\{\begin{array}{rl}
        h\cap F &;\,h^\ast\cap F\neq\varnothing\\
        F   &;\,h^\ast\cap F=\varnothing.
    \end{array}\right.
\end{equation}
We shall use the symbol $r$ to denote $r^X_F$ whenever the choices of $X$ and $F$ are unambiguous.
\begin{defn}\label{defn:restriction of a halfspace system to a subspace} Suppose $\HH$ is a halfspace system in the proper CAT(0) space $X$, and suppose $F$ is a closed convex subspace of $X$. We then define the restriction $\HH\big|_F$ of $\HH$ to $F$ as the set of all (possibly trivial) halfspaces of $F$ of the form $r^X_F(h)$ defined above for $h\in\HH$. 
\end{defn}
\begin{lemma}\label{lemma:restriction of a halfspace system} Suppose $\HH$ is a halfspace system in the proper CAT(0) space $X$, and $F$ is a closed convex subspace of $X$. Then the restriction $\HH\big|_F$ of $\HH$ to $F$ is indeed a halfspace system in the space $F$.
\end{lemma}
\proof{} Obviously, $\HH\big|_F$ forms a poc-set under inclusion and the restricted complementation operator. Let us denote $h\cap F$ by $h\big|_F$ for all $h\in\HH$, and let the image of any $S\subset\HH$ under the map $h\mapsto h\big|_F$ be denoted by $S\big|_F$. 

In order to make sure $\HH\big|_F$ is a halfspace system in $F$, we need to verify conditions $\mathbf{(H1),(H2)}$. We note that $S\big|_F$ is nested whenever $S\subset\HH$ is nested, and, accordingly, $S$ is transverse whenver $S\big|_F$ is transverse. Thus, applying axiom $\mathbf{(H2)}$ for $\HH$ we conclude there are no infinite transverse subsets in $\HH_F$, showing $\HH\big|_F$ is $\omega$-dimensional. The condition $\mathbf{(H1)}$ clearly holds.\ep\\

We now consider the dual $r^\circ$ of $r=r^X_F$. We note that $r$ is finite-to-one, except, possibly, for the preimages of $0,0^\ast\in\HH\big|_F$ being infinite: indeed, only a finite number of distinct walls of $\HH$ corresponds under $r$ to any given non-trivial wall of $\HH\big|_F$, because any family of halfspaces with intersecting walls in $\HH$ is transverse. This implies $r^\circ$ maps almost-equality classes of $\HH\big|_F^\circ$ into almost-equality classes of $\HH^\circ$, since for any pair $\sigma,\sigma'\in\HH\big|_F^\circ$, if $\sigma\vartriangle\sigma'$ is finite then 
\begin{equation} 
    r^\circ\sigma\vartriangle r^\circ\sigma'=r^{-1}(\sigma\vartriangle\sigma')
\end{equation}
is a finite set, too (since $0^\ast\in\sigma\cap\sigma'$, all the fibers of $r$ over $\sigma\vartriangle\sigma'$ are finite). Moreover, since $r$ is a surjection (by definition), we may conclude that if $r^\circ\sigma\vartriangle r^\circ\sigma'$ is finite, then $\sigma\vartriangle\sigma'$ must be finite.

Thus, $r^\circ$ induces an injective map $\Re\HH\big|_F\to\Re\HH$. This map will be of some importance to us, so let us fix some notation for it: 
\begin{defn} Suppose $\HH$ is a halfspace system in the proper CAT(0) space $X$, $F$ -- a closed convex subspace of $X$, and let $i_F:\bd F\to\bd X$ denote the inclusion map. Then the map $\Re\HH\big|_F\to\Re\HH$ induced by $r^\circ$ will be denoted by $\Re(i_F)$.
\end{defn}
We use the map $\Re(i_F)$ to understand the relation between the Roller decomposition maps associated with $F$ and $X$: 
\begin{prop}\label{prop:the restriction equality} Suppose $\HH$ is a halfspace system in the proper CAT(0) space $X$, and $F$ is a closed convex subspace of $X$. Let $\rho$ and $\rho_F$ denote the Roller decomposition maps associated with $X$ and $F$, respectively. Then, the following diagram is commutative:
\begin{equation}
\begin{CD}
    \bd F       @>\rho_F>>  \Re\HH\big|_F\\
    @Vi_FVV                @VV\Re(i_F)V\\
    \bd X       @>\rho>>    \Re\HH
\end{CD}
\end{equation}
\end{prop}
\proof{} We need to show $\Re(i_F)\circ\rho_F=\rho\circ i_F$. For this we fix a point $\xi\in\bd F$, and consider some $h\in\HH$: obviously, if $h\cap F$ contains $\xi$ in its boundary, then $h$ contains $\xi$ in its boundary; thus, if $\alpha\in\HH\big|_F^\xi$, then $r^\circ\alpha=r^{-1}(\alpha)$ lies in $\HH^\xi$, and we have 
\begin{equation}
    r^\circ(\HH\big|_F^\xi)\subseteq\HH^\xi,
\end{equation}
implying $\rho(\xi)\leq(\Re(i_F)\circ\rho_F)(\xi)$.

To prove equality, assume, by way of contradiction, that $\Re(i_F)(\rho_F(\xi))>\rho(\xi)$, and let $\alpha\in r^\circ(\rho_F(\xi))$. Since $\alpha\in\HH^\xi$, $\alpha\subset T(\xi)\cup P(\xi)$, and there exists an infinite descending chain $(c_n)_{n=1}^\infty$ of elements of $P(\xi)$ in $\alpha$ (as opposed to any element of $\rho(\xi)$, which may not contain such a chain). If $\alpha_n$ is a sequence of elements in $\rho(\xi)$ converging on $\alpha$, then we may assume $\alpha_n$ contains $c_n$ for all $n$.

Consider the chain $\bar c_n=r(c_n)$. Either $c_n\cap F$ is a dense (and open) subset of $F$ (possibly the whole of $F$) for all $n$, or $c_n\cap F$ and $c^\ast_n\cap F$ are both non-empty for $n$ large enough. In the first case, taking a point $p\in F$ and a point $q\in c_1^\ast$ we find out that the geodesic segment $[p,q]$ in $X$ intersects each of the walls $W(c_n)$, contradicting condition $\mathbf{(H1)}$ for $\HH$. Thus, only the second case is possible, but then, taking $\bar\alpha\in\rho_F(\xi)$ such that $\alpha=r^\circ(\bar\alpha)$, we find out that $(r(c_n))_{n=N}^\infty$ is a descending chain in $\bar\alpha$ for some $N$. This, too, is impossible, because no element of $\rho_F(\xi)$ may contain a descending chain of elements in $r(P(\xi))$, by the minimality of $\rho_F(\xi)$ (see proposition \ref{prop:boundary decomposition map is well defined}).\ep\\

\section{Connectivity properties of uniform halfspace systems.}\label{section:connectivity properties} From now on assume $G$ is a group acting properly-discontinuously and co-compactly by isometries on $X$, keeping a uniform (and $\omega$-dimensional!!) halfspace system $\HH$ invariant.

\subsection{Lines and flat sectors.}
We recall a standard fact about CAT(0) spaces satisfying our assumptions:
\begin{prop}[\cite{[BH]} II proposition 9.5(3)]\label{prop:angles are realized} For all $\xi,\xi'\in\bd X$ there exists a point $p\in X$ and representative rays $c\in\xi$ and $c'\in\xi'$ emanating from $x_0$ such that $\angle(\xi,\xi')=\angle_p(c,c')$. In particular, if $\angle(\xi,\xi')=\pi$, then there exists a geodesic line $c:\RR\to X$ such that $\xi=c(\infty)$ and $\xi'=c(-\infty)$. 
\end{prop}
\begin{lemma}\label{lemma:rho separates the ends of a line} Suppose $\rho^{-1}(\Pi)=\varnothing$, and suppose $c:\RR\to X$ is a geodesic line with endpoints $\xi=c(\infty)$ and $-\xi=c(-\infty)$. Then $\rho(\xi)$ and $\rho(-\xi)$ are incomparable.
\end{lemma}
\proof{} Given the line $c$, we will show that $\rho(\xi)\leq\rho(-\xi)$ implies $\rho(\xi)=\Pi$, producing a contradiction.

For any $h\in T(\xi)$, note that if $c(\RR)\nsubseteq\cl{h}$, then $-\xi$ cannot lie in $\bd W(h)$, implying $h^\ast\in T(-\xi)$, which is the same as $h\in T(-\xi)^\ast$ and implies $h\notin T(-\xi)$.

Now, the idea of the proof is to use the fact that $T(\xi)$ is non-empty to produce an element $h\in T(\xi)\minus T(-\xi)$, making the inequality $\rho(\xi)\leq\rho(-\xi)$ impossible: if this inequality holds, then corollary \ref{cor:comparison in bdd type} implies $T(\xi)$ is contained in $T(-\xi)$, contradicting the existence of such an $h$. 

In view of the remark above we may assume that every $h\in T(\xi)$ closure-contains the line  $c(\RR)$. We may also assume no element of $T(\xi)$ lies in $P(-\xi)$, because this would also end the proof. Thus, we continue operating under the assumption that $c(\RR)$ is contained in $h$, for all $h\in T(\xi)$. But this is impossible: by uniformness, since $T(\xi)$ is non-empty there exists an infinite descending chain $(c_n)_{n=1}^\infty$ in $T(\xi)$; picking a point $p\in c_1^\ast$ and a point $q\in c(\RR)$, we immediately obtain a contradiction to $\mathbf{(H1)}$.\ep
\begin{cor}[$\rho$ separates $\pi$-discrete sets]\label{cor:rho separates pi-discrete sets} Suppose $\rho^{-1}(\Pi)=\varnothing$, and $A$ is a subset of $\bd X$ which is $\pi$-discrete in the angular metric. Then the restriction of $\rho$ to $A$ is injective.\ep
\end{cor}
The last corollary shows that the Roller boundary of a uniform halfspace system is particularly interesting when $\bd X$ is, say, connected with respect to the angular metric. 
\begin{example} Suppose $X=\FAT{H}^n$, and fix $x_0\in X$. For any $x\in X$ and $g\in G$ let  $h_{x,g}$ denote the halfspace of all points $p\in X$ which are closer to $x$ than to $gx$. Then the set of halfspaces $\HH$ consisting of all $h_{x,g}$ with $x\in G\cdot x_0$ and $g\in G$ is a uniform halfspace system in $X$ satisfying $\rho^{-1}(\Pi)=\varnothing$. Since any two points of $\bd X$ are a distance $\pi$ apart in the angular metric, we conclude $\rho$ is one-to-one. 

Thus, in order for $\rho$ to be meaningful for a group $G$ acting on $\FAT{H}^n$, one needs $G$ to have parabolic points (and then $\HH$ stops being uniform). Note that the same kind of problem will arise for any visibility space, whenever $\HH$ is rich enough to satisfy $\rho^{-1}(\Pi)=\varnothing$.
\end{example}

\subsection{Improved closure formula.}
Let us now revisit the formula $\mathbf{(FF_0)}$ we had derived in the previous section. It turns out that now we are able to make it into an equality.
\begin{thm}\label{thm:closure formula 1} Suppose $G$ is a group acting properly and co-compactly by isometries on a proper CAT(0)-space $X$, and $\HH$ is a $G$-invariant uniform halfspace system in $X$. Then, for all $\Sigma\in\rho(\bd X)$ one has the equality 
\begin{equation}
    \mathbf{(FF)}\quad\cl{\rho^{-1}(\Sigma)}=\bigcup_{\Sigma_1\leq\Sigma}\rho^{-1}(\Sigma_1)\,.
\end{equation}
\end{thm}
\proof{} Suppose $\Sigma_1<\Sigma$, and consider a boundary point $\xi$ satisfying $\rho(\xi)=\Sigma_1$. 

Let us write $\Sigma=\rho(\eta)$ for some $\eta\in\bd X$. By corollary \ref{cor:comparison in bdd type} and \ref{cor:comparison in bdd type - 2}, we have 
\begin{enumerate}
    \item $P(\eta)\subset P(\xi)$, and --
    \item $T(\xi)\subset T(\eta)\subset T(\xi)\cup P(\xi)$.
\end{enumerate}
In particular, $T(\eta)$ is non empty.

Now, since $\rho(\eta)$ and $\rho(\xi)$ are comparable, lemma \ref{lemma:rho separates the ends of a line} implies we have $\angle(\xi,\eta)<\pi$. Since there exists in $X$ a point $p$ such that $\angle_p(\xi,\eta)=\angle(\xi,\eta)$, we conclude that the ideal triangle $\vartriangle(p,\xi,\eta)$ in $X$ bounds a flat sector $F$ in $X$ (see \cite{[BH]}, cor.9.9, p.283). For any boundary point $\zeta$, let us denote the geodesic ray from $p$ to $\zeta$ by $[p,\zeta]$. 

Let us now consider any point $\zeta\in\bd F$ other than $\xi$. If we show that $\rho(\zeta)=\rho(\eta)$, then we are done, since $\xi$ may clearly be written as the limit of a sequence of such $\zeta$.

We fix $\zeta\in\bd F$ and note the ray $[p,\zeta]$ lies in the convex hull $F$ of $[p,\xi]$ and $[p,\eta]$. 

First, for any $h\in P(\eta)=P(\eta)\cap P(\xi)$ we note that since there exists $A>0$ such that both $[p,\xi]$ and $[p,\eta]$ are eventually-contained in $N_A(W(h))$, so does $[p,\zeta]$, by the convexity of the metric in $X$; this proves $P(\eta)\subseteq P(\zeta)$. Next, if $h\in T(\eta)$, then from $h\in T(\xi)\cup P(\xi)$ we deduce in a similar fashion that $h$ cannot lie in $T(\zeta)^\ast$, and we have shown $T(\eta)\subset T(\zeta)\cup P(\zeta)$, proving $\rho(\zeta)\leq\rho(\eta)$.\\

Let us now use $\rho(\xi)\leq\rho(\eta)$ more explicitly: there necessarily exists an infinite descending chain $(c_n)_{n=1}^\infty$ of elements in $P(\xi)\cap T(\eta)$. Selecting points $p_n\in[p,\eta]\cap c_n$ we necessarily have that $(p_n)_{n=1}^\infty$ converges on $\eta$ -- otherwise there would be a geodesic segment $[p,q]\subset[p,\eta]\cap X$ intersecting all the walls $W(c_n)$. Letting $\gamma_n=[p_n,\xi]$, we observe that $[p,\zeta]$ intersects each $\gamma_n$ in a point $q_n\in F$, and that these points converge on $\zeta$. 

Now consider $h\in T(\xi)$. We already have $T(\xi)\subset T(\eta)\subseteq T(\zeta)\cup P(\zeta)$. Assume $h\in P(\zeta)$. Then, there exists $B>0$ such that all but finitely many points of the sequence $q_n$ belong to $N_B(h^\ast)$. By the convexity of the metric, all the rays $\gamma_n$ are then eventually contained in $N_B(h^\ast)$, contradicting $h\in T(\xi)$. Thus we have shown $T(\xi)\subseteq T(\zeta)$ holds, which also implies $T(\zeta)\subseteq T(\xi)\cup P(\xi)$ and $P(\zeta)\subset P(\xi)$.\\

Finally we are in position to show $T(\eta)\subseteq T(\zeta)$, which will finish the proof. We already have $T(\eta)\subseteq T(\zeta)\cup P(\zeta)$, so let us assume there exists $T(\eta)\cap P(\zeta)\neq\varnothing$ and find out what goes wrong. We use the equality $P(\zeta)=T(\rho(\zeta))$ to deduce there exists an infinite descending chain of such elements in $T(\eta)\cap P(\zeta)$ -- call it $(h_n)_{n=1}^\infty$. Without loss of generality, $F\not\subseteq h_1$: if $F$ were contained in all the $h_n$, then taking a point $p'\in h_1^\ast$, the geodesic segment $[p,p']$ would have intersected $W(h_n)$ for all $n$, which is impossible.

Since $h_n\in P(\zeta)\subseteq P(\xi)$, the intersection $F_n=F\cap\cl{h_n}$ is itself a flat sector in $X$ having the same ideal boundary as $F$. Consider the ray $[p,\zeta]$: for large enough $n$ we have $p\in h_n^\ast$ and $[p,\zeta]$ intersects $W(h_n)$ transversely; by the convexity of the metric, $h_n\in P(\zeta)$ is impossible, and we are done.\ep

\subsection{Boundary paths}\label{subsection:boundary paths} Inspecting the above proof one immediately notes the connection between the structure of the Roller boundary of $\HH$ and properties of ``Euclidean'' paths in the boundary. Adapting an idea of Croke and Kleiner (\cite{[CroKle]}), we define a special type of paths in the boundary of $X$ as follows:
\begin{defn}[safe paths] A path $\alpha:[0,1]\to\bd X$ is said to be \emph{safe} with respect to $\HH$, if $\rho\circ\alpha([0,1])$ is a finite subset of $\Re\HH$.
\end{defn}
Obviously, the concatenation of two safe paths is again safe, so that safe path-components of $\bd X$ are defined (with respect to a particular choice of $\HH$). 
\begin{lemma} safe path components of $\bd X$ are saturated with respect to $\rho$, i.e., for every safe path component $C$ of $\bd X$ we have $\rho^{-1}(\rho(C))=C$. 
\end{lemma}
\proof{} It suffices to prove that for any $\xi,\eta\in\bd X$, if $\rho(\xi)=\rho(\eta)$, then $\xi$ and $\eta$ may be joined by a safe path. The procedure described in the proof of the theorem above gives the required path. \ep
\begin{prop}[characterizing safe components]\label{prop:safe paths} Two points $\xi,\eta\in\bd X$ lie in the same safe path component if and only if there exists a sequence $\Sigma_0,\ldots,\Sigma_n$ of elements in $\rho(\bd X)$ satisfying
\begin{enumerate}
    \item $\rho(\xi)=\Sigma_0$ and $\rho(\eta)=\Sigma_n$;
    \item $\Sigma_{i-1}$ and $\Sigma_i$ are comparable for all $i\in\{1,\ldots,n\}$.
\end{enumerate}
We call such a sequence a \emph{connecting sequence of boundary classes of length $n$} for the points $\xi$ and $\eta$. 
\end{prop}
\proof{} Once again, given the points $\xi,\eta$, the procedure used for the proof of $\mathbf{(FF)}$ shows that the existence of such a ``connecting sequence'' in $\rho(\bd X)$ implies there is a safe boundary path from $\xi$ to $\eta$. For the converse, we claim that for any path $\alpha:[0,1]\to\bd X$, the set $\rho(\mathrm{Im}\alpha)$ -- we shall abbreviate it as $\rho(\alpha)$ -- contains a sequence of boundary classes connecting its endpoints.

Suppose not -- then there exists a path providing a counter-example to our claim with $\rho(\alpha)$ of minimal possible size. 

Let $A$ be the set of all $i\in\{1,\ldots,n\}$ such that $\Sigma_j\leq\Sigma_i$ holds for no $j\neq i$. Then $\mathbf{(FF_0)}$ and the continuity of $\alpha$ imply the set $T_i=\alpha^{-1}(\rho^{-1}\Sigma_i)$ is a closed subset of $[0,1]$. Let $J_i$ denote the closed subinterval of $[0,1]$ spanned by $T_i$. 

Among all $i\in A$ find those with $J_i$ maximal with respect to inclusion. on each such interval we use the previous lemma to redefine $\alpha$ in such a way that $\rho(\alpha(t))=\Sigma_i$ for all $t\in J_i$. Thus, by the minimality property of $\alpha$, all the $J_i$ are disjoint. 

For each $i\in A$, set $t(i)=\min(J_i)$, $t'(i)=\max(J_i)$, and let us write $A\minus\{0\}=\{i_1,\ldots,i_k\}$ with $t(i_1)<\ldots<t(i_k)$. Thus, we may break the interval $[0,1]$ into a series of consecutive intervals 
\begin{equation}
    \begin{array}{rclrcl}
        I_0&=&[0,t(i_1)],& I'_1&=&[t(i_1),t'(i_1)]=J_{i_1}\\
        I_1&=&[t'(i_1),t(i_2)],& I'_2&=&[t(i_2),t'(i_2)]=J_{i_2}\\
        &\vdots\\
        I_k&=&[t'(i_k),1],
    \end{array}
\end{equation}
such that each interval $I'(\ell)$ is mapped by $\alpha$ into $\rho^{-1}(\Sigma_{i_\ell})$. 

Suppose $k>1$. Then for each of the intervals $I_\ell$ ($0\leq\ell\leq k$) the restriction $\alpha_\ell$ of $\alpha$ to $I_\ell$ is a path satisfying $\rho(\alpha_\ell)$ is a proper subset of $\rho(\alpha)$, showing $\rho(\alpha_\ell)$ must contain a sequence of boundary classes connecting its endpoints. But then $\rho(\alpha)$ contains a sequence of boundary classes containing $\alpha(0)$ and $\alpha(1)$, which is impossible. Thus we deduce $k\geq 1$, implying $\rho(\alpha)$ has at most one minimal element -- denote it by $\Sigma_1$ -- except $\Sigma_0$ (which may, or may not be minimal).

Note now that $\Sigma_0\leq\Sigma_n$ is impossible, by our assumption on $\alpha$. Thus, $\Sigma_1$ does indeed exist ($k$ cannot be zero).

Suppose now that $\Sigma_1$ is defined and $\Sigma_0$ is a minimal element of $\rho(\alpha)$. Then $[0,1]=I_0\cup I'_1\cup I_1$ and $\rho(\alpha)$ may be written as the union of the two sets
\begin{equation}
    B_0=\left\{\Sigma\in\rho(\alpha)\,\big|\,\Sigma\geq\Sigma_0\right\},\quad
    B_1=\left\{\Sigma\in\rho(\alpha)\,\big|\,\Sigma\geq\Sigma_1\right\},
\end{equation}
and assume $\Sigma_n\notin B_0$ (else we are already done).

Now, $\mathbf{(FF_0)}$ implies the sets $B_0, B_1$ are not disjoint, for otherwise $\rho^{-1}(B_0)$ and $\rho^{-1}(B_m)$ split $\mathrm{Im}\alpha$ as a union of two closed, non-empty disjoint sets, which is impossible.

Taking $\Sigma\in B_0\cap B_1$, consider the path $\alpha'=\alpha\big|_{I_0\cup I'_1}$: $\rho(\alpha')$ evidently contains the sequence $(\Sigma_0,\Sigma,\Sigma_1)$, which connects $\alpha(0)$ with $\alpha(t'(1))$. Next, the path $\alpha''=\alpha\big|_{I_1}$ has $\rho(\alpha'')$ contain the sequence $(\Sigma_1,\Sigma_n)$, which connects $\alpha(t'(1))$ with $\alpha(1)$. Thus, $\rho(\alpha)$ contains a sequence of boundary classes containing its endpoints -- a contradiction.

Thus, we have shown $\rho(\alpha)$ has a unique minimal element $\Sigma_1\neq\Sigma_0,\Sigma_n$, but then $(\Sigma_0,\Sigma_1,\Sigma_n)$ is a sequence connecting the endpoints of $\alpha$, implying our $\alpha$ simply does not exist.\ep\\

Thus, the safe path components of $\bd X$ correspond to comparability components in $\rho(\bd X)$, and the natural question to ask now is how many such comparability components are there in $\Re\HH$.

\subsubsection{Example: Flats.}\label{special case:full flats} Let $X$ be as before, and suppose $F$ is a flat -- i.e., an isometric copy of $\EE^m$ in $X$, for some positive $m$. We claim that $\bd F$ is safely path-connected. By lemma \ref{lemma:restriction of a halfspace system}, $\HH\big|_F$ is a halfspace system on $F$, and \ref{prop:the restriction equality} tells us that $\rho_F(\bd F)$ is in one-to-one correspondence with the image of $\bd F$ under $\rho$ in $\Re\HH$. Thus it will be enough to show that $\bd\EE^m$ is safely path-connected with respect to any halfspace system on $\EE^m$. Next, since geodesics in $\bd\EE^m=\SS^{m-1}$ are arcs of great circles, another application of the same two lemmas implies it is enough to show our claim for the case $m=2$. 

Suppose $\HH$ is a halfspace system on $\EE^2$, and for each $\alpha\in[0,\pi)$ denote by $\HH_\alpha$ the subset of those $h\in\HH$ such that the (positive) angle from the $x$-axis to $W(h)$ equals $\alpha$. In $\EE^2$ we necessarily have that, for all $\alpha,\beta\in[0,\pi)$ and all $h\in\HH_\alpha,k\in\HH_\beta$, $h\pitchfork k$ holds if and only if $\alpha\neq\beta$; since $\HH$ contains no infinite transverse subset, we conclude $\HH_\alpha$ is empty for all but finitely many values of $\alpha\in[0,\pi)$. Denote those values by $\{\alpha_1,\ldots,\alpha_d\}$, ordered in increasing order. With each $\alpha_i$ we associate two boundary points -- denote them by $\pm\xi_i$, -- the two endpoints of a line intersecting the $x$-axis at a (positive) angle $\alpha_i$. 

Now, for any point $\xi\in\bd\EE^2$, if $\ell$ is a line  with endpoint $\xi$, then for any $\alpha\in[0,\pi)$ and any $h\in\HH_\alpha$ we have $h\in P(\xi)$ if and only if the angle of $\ell$ to the $x$-axis (measured from the positive ray of the $x$-axis to $\ell$) equals $\alpha$. Thus, $P(\xi)$ is non-trivial if and only if $\xi\in\{\pm\xi_i\}_{i=1,\ldots,d}$. Moreover, as $\xi$ moves along the interior of a boundary arc whose endpoints are a pair of consecutive points of the cyclically-ordered set $\{\pm\xi_i\}_{i=1,\ldots,d}$, $T(\xi)$ remains constant, showing that $\rho(\xi)$ also remains constant. Thus, the image of $\bd\EE^2=\SS^1$ under $\rho$ consists of at most $4d$ elements, showing $\bd\EE^2$ is safely path-connected. We have proved: 
\begin{prop}[Euclidean boundaries are safe]\label{prop:Euclidean boundaries are safe} Suppose $X$ is a proper CAT(0) space and $\HH$ is a halfspace system on $X$. Then, for any isometrically embedded flat $F\cong\EE^m$ in $X$, its boundary $\bd F$ in $\bd X$ is safely-connected with respect to $\HH$, with $\rho(\bd F)$ isomorphic to $\Re(\EE_{\square}^M)$ for some $M\geq m$, where $\EE^M$ is taken with the standard halfspace system.\ep 
\end{prop}
\begin{remark} We have just proved that the Roller boundary of a halfspace system in a two-dimensional flat is finite. Using normal vectors (rather than angles) it is possible to generalize the same argument to prove that the Roller boundary of a halfspace system in a flat of any finite dimension is finite.
\end{remark}
In the case when $X$ and $\HH$ are invariant under a proper co-compact action by a group $G$ of isometries of $X$, one way of obtaining a flat plane $F$ is to find a subgroup $A$ of $G$ isomorphic to $\ZZ^2$ and then applying the Flat Torus Theorem to produce a flat $F$ which is invariant under $A$. In this case one easily sees that the invariance of $\HH$ under $A$ implies the invariance of all the $\HH_{\alpha_i}$ under $A$, showing $\HH$ is poc-set isomorphic to the standard cubulation of $\EE^d$. It is easy to see that the Roller boundary of this cubulation is isomorphic to the barycentric subdivision of a $d$-dimensional cube, and the image of $\bd F$ under $\rho$ then traverses the cycle of length $4d$ spanned by the classes $\rho(\pm\xi_i)$.

\subsubsection{Example: Flat sectors and uniform systems.} The preceding example explains the role of flats for a general halfspaces system, showing that boundaries of flats are always safe. That argument does not work for flat {\it sectors}, however, but it can be mended under the additional assumption that $\HH$ satisfies the parallel walls property -- in particular, when $\HH$ is a uniform system invariant under a geometric group action.

Observe that for every $h\in\HH$ for which the restriction to $F$ is proper, the wall $W(h)$ intersects $F$ in an interval separating $F$. The tool that makes the parallel walls property relevant is the following:
\begin{lemma} Suppose $F$ is a flat sector and $h\in\HH$ restricts to a proper halfspace of $F$. Suppose that $\gamma$ is a geodesic ray in $F$ such that $d(\gamma(t),F\cap\cl{h^\ast})$ is unbounded. Then $d(\gamma(t),h^\ast)$ is unbounded.
\end{lemma}
\proof{} Let $v$ be the vertex point of $F$. If $d(\gamma(t),h^\ast)$ were bounded, that would imply $h\in P(\gamma(\infty))$. Let $x\in W(h)\cap F$ be any point. Since $W(h)$ is convex and complete, we must have that $[x,\gamma(\infty))$ is contained in $W(h)$. However, $x\in F$ and $\gamma\subset F$; since $F$ is convex and complete, this implies $[x,\gamma(\infty))$ lies in $F$. Thus, $[x,\gamma(\infty))$ lies in $W(h)\cap F$, contradicting the assumption regarding $\gamma$.\ep
\begin{prop} Suppose $X$ is a proper $CAT(0)$ space and $\HH$ is a halfspace system in $X$, satisfying the parallel walls property. If $F$ is a flat sector in $X$, then $\bd F$ is safely-connected with respect to $\HH$.
\end{prop}
\proof{} It is enough to prove the proposition for the case when $F$ has a vertex angle $0<\theta<\pi$. We will identify $F$ with the set 
\begin{equation*}
	\left\{re^{i\alpha}\in\CC\,\big|\,r\geq 0\,,\quad 0\leq \alpha\leq\theta\right\}\,,
\end{equation*}
with $0$ corresponding to the vertex of $F$ in $X$. To each angle in the interval $[0,\theta]$ corresponds a unique point $\xi_\alpha\in\bd F$. Same as before, we decompose (the proper part of) $\HH\big|_F$ as the union of subsystems $\HH_\alpha$ (throw in the trivial halfspaces) where $\alpha\in[0,\pi]$ and $h\in\HH_\alpha$ if and only if the angle between the interval $W(h)$ and the $x$-axis (of $F$) equals $\alpha$. Then it is clear that $\HH_\beta\subseteq T(\xi_\alpha)$ whenever $\beta\neq\alpha$, and that $H_\alpha=P(\xi_\alpha)$. Thus, in order for $\bd F$ to be safely connected it is sufficient that $\HH_\alpha$ be trivial for all but finitely many values of $\alpha\in(0,\theta)$.

For every $\alpha\in(0,\pi)$ and $h\in\HH_\alpha$ we will denote the point of intersection of $W(h)$ with the boundary rays of $F$ by $h(0)$, and the endpoint of the ray arising as $W(h)\cap F$ will be denoted by $h(\infty)$.

The assumption that $\HH$ has no infinite transverse subset is used as follows. Suppose $h_n\in\HH_{\alpha_n}$ ($n\in\NN$) satisfy
\begin{enumerate}
	\item $h_n(0)\in\RR$ for all $n$,
	\item $h_{n+1}(0)>h_n(0)$ for all $n$,
	\item $\alpha_{n+1}>\alpha_n$ for all $n$.
\end{enumerate}
Then it is clear that the rays $W(h_n)\cap F$ cross pairwise. Thus, such a configuration is impossible. A symmetric configuration (for which the rays $W(h_n)\cap F$ have $h_n(0)=r_ne^{i\theta}$, and with the $\alpha_n$ {\it decreasing}) is, of course, also impossible for the same reasons.

Consider the set $A$ of all $\alpha\in(0,\theta)$ for which $\HH_\alpha$ is poc-isomorphic to the standard halfspace system on $\EE^1$. If $A$ is infinite, extract a strictly monotone sequence $\alpha_n$ from $A$; because of the discreteness assumption on $\HH$, it is then easy to construct (inductively) a forbidden configuration in $\HH\big|_F$, producing a contradiction.

Therefore, it will be enough to prove that $\HH_\alpha$ is poc-isomorphic to the standard halfspace system on $\EE^1$ whenever it is non-trivial. From symmetry considerations it will be enough to prove that every $a\in\HH_\alpha$ eventually-containing the ray of (positive) reals has a $b\in\HH_\alpha$ satisfying $b<a$.

Fix $a$ as above, and consider the ray of points of $F$ we had identified with $\RR_+$. By the preceding lemma, a point on this ray sufficiently far away from $W(a)\cap F$ must be contained in an element $h_1\in\HH$ satisfying $h_1<a$. This implies that $h_1$ restricts to a proper halfspace of $\HH\big|_F$, so that $h_1\in\HH_{\alpha'_1}$ for $\alpha'_1\leq\alpha$. We may take $h_1$ to be maximal with this property (as the interval $[h_1,a]$ is finite). In particular, there are no intermediate halfspaces between $a$ and $h_1$.

If $\alpha'_1=\alpha$, we are done. If not, Then we apply the preceding lemma again, to the halfspace $a$ and the ray $W(h_1)\cap F$. This results in a halfspace $h_2<a$, $h_2\in\HH_{\alpha_2}$ and since $a$ and $h_1$ had no intermediate halfspaces, we also obtain $\alpha_2>\alpha_1$ and $h_2(0)>h_1(0)$. Proceeding inductively in the same manner, we obtain a forbidden configuration, as desired.\ep\\

As a corollary, we obtain the following result:
\begin{thm}[Tits components and safe components coincide]\label{thm:Tits components computed} Suppose $G$ is a group acting geometrically on a CAT(0) space $X$, and suppose $\HH$ is a $G$-invariant uniform halfspace system. Then, for every $\xi\in\bd X$, the open Tits ball of radius $\pi$ about $\xi$ is contained in the safe component of $\xi$. In particular, the components of the Tits boundary of $X$ coincide with the safe components of $\bd X$.\ep
\end{thm}
Thus, the decomposition map $\rho$ defined by $\HH$ introduces a new structure on the Tits boundary, determined by the ordering of the Roller boundary associated with $\HH$.

\subsubsection{Example: $\FAT{H}^2\times\FAT{H}^2$} Let $X_1,X_2$ be visible proper CAT(0) spaces, and let $X=X_1\times X_2$. We claim: 
\begin{prop} If $\HH$ is \emph{any} uniform halfspace system on $X$, then $\bd X$ is safely-connected with respect to $\HH$.
\end{prop}
\proof{} the boundary $\bd X$ is naturally homeomorphic to the spherical join $\bd X_1\ast\bd X_2$. Let $\xi=[\xi_1,\xi_2,\alpha]$ and $\eta=[\eta_1,\eta_2,\beta]$ be two distinct points of $\bd X$ (with their standard representations as points of $\bd X_1\ast\bd X_2$, where $\xi_i,\eta_i\in\bd X_i$ and $\alpha,\beta\in[0,\pi/2]$).

Now, there are geodesic lines $\ell_i\subset X_i$ joining $\xi_i$ to $\eta_i$, and we may consider the embedded flat $F=\ell_1\times\ell_2\subset X$: we note that $\bd F$ is then the join of the two-point subspaces $\{\xi_1,\eta_1\}$ and $\{\xi_2,\eta_2\}$. In particular, $\bd F$ contains our points $\xi$ and $\eta$. By proposition \ref{prop:Euclidean boundaries are safe}, there is a safe path in $\bd X$ from $\xi$ to $\eta$.\ep\\ 

We would like to draw the attention of the reader to the fact that, in this example, the Roller boundary of $\HH$ may be considered as an extremely highly-connected graph when both boundaries $\bd X_1$ and $\bd X_2$ are infinite. For example, given any finite set of points $A$ in $\bd X$, and any pair of points $\xi,\eta\in\bd X\minus A$, there will be a safe path from $\xi$ to $\eta$ missing $A$, and this is independent of the choice of $\HH$ (so long as $\HH$ is uniform). On the other hand, $\Re\HH$ is by no means a trivial structure: If $X$ and $\HH$ admit a geometric action by some group, then since the equatorial copies of the $\bd X_i$ in $\bd X$ are infinite $\pi$-discrete sets, $\rho$ is injective on each of them, producing pairwise-incomparable classes; this shows that, as a graph, $\rho(\bd X)$ may be quite complicated (an infinite complete bipartite graph in this case).

\subsubsection{Example: the Croke-Kleiner examples}\label{subsubsection:Croke-Kleiner} We recall the set of examples by Croke and Kleiner (\cite{[CroKle]}).

Given $\alpha\in(0,\pi/2]$, let $T_1,T_2$ be two standard (``square'') geometric $2$-tori (with the standard CW decomposition), and let $T_0$ be the geometric $2$-torus obtained from a flat rhombus $R$ with an interior angle $\alpha$ via the standard gluing; the two simple closed curves arising as the image of $\partial R$ in $T_0$ will be denoted by $m'_1$ and $m'_2$. Selecting meridional curves $m_1,m_2$ of unit length in $T_1,T_2$ respectively, we glue each $T_i$ ($i\in\{1,2\}$) to $T_0$ by identifying $m_i$ with $m'_i$ using an isometry.  Let now $X_\alpha$ denote the universal cover of the resulting ``torus complex'' $Y_\alpha$, and use Bridson's theorem to metrize $X$ as a piecewise-Euclidean cell complex. It is known that $X_\alpha$ is a CAT(0) space, and we have that $\pi_1(Y_\alpha)$ which is independent of the  choice of $\alpha$) acts properly discontinuously and co-compactly by isometries on $X_\alpha$. We let $\HH$ be the halfspace system arising naturally from the cube structure on $X_\alpha$, and it is evident that $\HH$ is uniform, and suppress the index $\alpha$ until it is needed. 

In their paper \cite{[CroKle]}, Croke and Kleiner define a system $\hsm{W}$ of \emph{walls} in $X$ -- the distinct (disjoint) lifts of $T_0$ --, and a system $\hsm{B}$ of \emph{blocks} -- all lifts of subsets of $Y$ of the form $T_i\cup_{m_i=m'_i}T_0$. Let us color the $1$-skeleton of $X$ as follows: edges of $X$ projecting to the curves $m_i=m_i'$ will be colored by $0$, while edges projecting to curves not contained in $T_0$ will be colored according to their projections being contained in $T_1$ or $T_2$. The $1$-skeleton of each block $B$ is then colored using two colors (one of which is necessarily $0$), according to the way it was constructed.

\begin{figure}[ht]
    \centering\includegraphics[width=4 in]{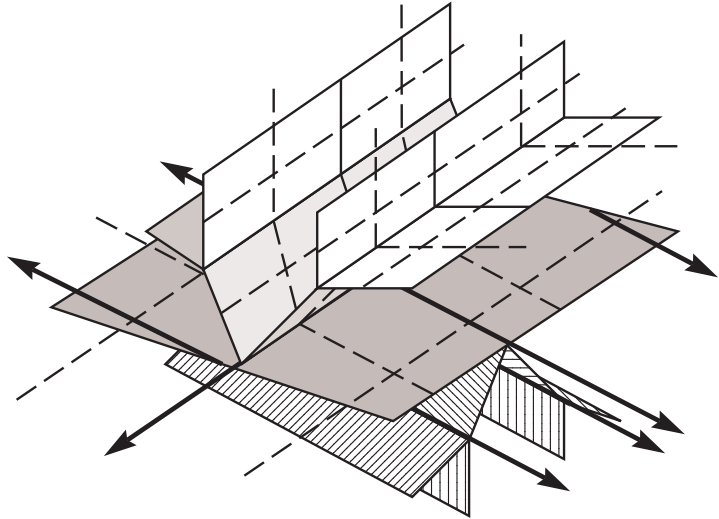}
    \caption{\protect\scriptsize Intersection patterns of halfspaces and blocks in the
Croke-Kleiner examples: the horizontal plane is the intersection
plane of two neighbouring blocks (above and below the plane,
respectively); singular fibers are denoted by arrows; walls (of
the cubing) are denoted by dashed lines.}
    \label{figure:calculating halfspaces in Cr-Kl}
\end{figure}

We shall now list some facts from \cite{[CroKle]}. 
\begin{description}
    \item[local properties of blocks.] Each block $B$ is isometric to the cartesian product of a $4$-regular metric tree -- all of whose edges have length $1$ -- with $\RR$, which is also standardly realized as a metric tree, whose vertex set equals $\ZZ$. 
    
    The fibers of the projection to the tree factor are called \emph{singular fibers}. Note they are all geodesic lines in $X$, all parallel to each other.
    \item[global properties of blocks.] The interaction among distinct blocks is as follows:
    \begin{enumerate}
        \item $\hsm{B}$ covers $X$;
        \item Two distinct blocks are either disjoint or share precisely one wall, which separates them in $X$;
        \item The nerve of the covering $\hsm{B}$ of $X$ is a tree (actually, the Bass-Serre tree of the splitting of $\pi_1(X)$ as an amalgam over $\pi_1(T_0)$). Put more simply, the graph whose vertex set is $\hsm{B}$ and whose edge set is $\hsm{W}$, with an edge $W$ incident to an edge $B$ iff $W\subset B$, is a tree. We shall need the separation order on this tree, and write \sep{a}{b}{c} to denote that an element $b$ of this tree separates the element $a$ from the element $c$.
    \end{enumerate}
\end{description}
Using these facts, we may now proceed to give a precise description of $\HH$. Given the structure of an individual block $B$, if a halfspace $h\in\HH$ is such that both $h$ and $h^\ast$ intersect $B$, then either $h\cap B$ has the form of the cartesian product of a halfspace of $\RR$ with the $4$-regular tree corresponding to $B$, or $h\cap B$ is the product of a halfspace in the $4$-regular tree with the whole of $\RR$. In the former case, let us write $h\in T(B)$, while in the latter we shall write $h\in P(B)$, according to $h$ being transverse or parallel to the singular fiber of the block $B$. $P(B)$, in turn, splits as the disjoint union of two sets $P_c(B)$ and $P_0(B)$, depending on whether $W(h)$ projects (under the covering $X\to Y$) into a non-zero (colored) edge, or into a zero-colored edge.  Suppose now that $B$ is a block. We consider the three cases discussed above for $h\in\HH$ with respect to $B$, assuming neither $h$ nor $h^\ast$ contain $B$. The discussion is based on the intersection patterns among the liftings of the tori $T_1,T_2$ and $T_0$, as pictured in figure \ref{figure:calculating halfspaces in Cr-Kl}.

\begin{itemize}
    \item $h\in T(B)$.\\ In this case, for any block $B'$ adjacent to $B$, we must have $h\in P_0(B')$, and $W=B\cap B'$ is the only wall of $B'$ intersecting $W(h)$, because the intersection of $W(h)$ with $W$ is a singular fiber of $B'$. In particular, for every $B''\in\hsm{B}$ satisfying \sep{B}{B'}{B''} we must either have $B''\subset h$ or $B''\subset h^\ast$. 
    
    Thus, $T(B)$ induces non-trivial halfspace systems only on $B$ and its immediate neighbours. 
    \item $h\in P_c(B)$.
    
    In this case, $W(h)$ is equal to a singular fiber of $B$. If $B'\neq B$ is any block and $W$ is the unique wall of $B$ satisfying \sep{B}{W}{B'}, then we must have either $B'\subset h$ or $B'\subset h^\ast$ according to which of $W\subset h, W\subset h^\ast$ (resp.) occurs.
    \item $h\in P_0(B)$.
    
    This case is actually symmetric to the first one: there is a unique wall $W$ of $B$ containing the singular fiber $W(h)\cap B$, and therefore $h\in T(B')$, where $B'$ is the unique block satisfying $B\cap B'=W$. 
\end{itemize}
From this analysis we immediately deduce that to every proper element $h\in\HH$ there corresponds a unique block $B\in\hsm{B}$ satisfying either $h\in T(B)$ or $h\in P_c(B)$. Let us denote this block by $B(h)$, and decompose $\HH$ as the disjoint union of subfamilies \begin{equation}
    \HH_B=\left\{h\in\HH\,\big|\,B(h)=B\right\}\cup\{\varnothing,X\}.
\end{equation}
We note that $\HH_B\cap\HH_{B'}$, where $B,B'$ is a pair of adjacent blocks, produces a halfspace system $\HH_W$ (corresponding to the wall $W=B\cap B'$) such that the restriction map $r^X_W:\HH_W\to\HH\big|_W$ is an isomorphism, since for every $h\notin\HH_W$ we have either $W\subset h$ or $W\subset h^\ast$.

We need now to compute $\HH^\circ$ and the map $\rho$. Suppose $\alpha\in\HH^\circ$ is not principal, so that $\alpha$ contains an infinite descending chain $(h_n)_{n=1}^\infty$; denote $B_n=B(h_n)$, and consider the following cases: 
\begin{description}
	\item[The sequence $(B_n)_{n=1}^\infty$ does not stabilize.] Passing to an (equivalent) subsequence we may assume that $B_n$ and $B_{n+1}$ are disjoint for all $n$. From the analysis above it follows that we may find, for each $n$, a block $B'_n$ satisfying $\sep{B_n}{B'_n}{B_{n+1}}$ and a halfspace $h'_n\in P_c(B'_n)$ satisfying $h_{n+1}<h'_n<h_n$.
     
	The sequence $(h'_n)_{n=1}^\infty$, apart from being equivalent to the original sequence $(h_n)_{n=1}^\infty$, also has the property that every $h\in HH$ must satisfy either $h'_n<h$ for sufficiently large $n$, or $h'_n<h^\ast$ for sufficiently large $n$.  Consequently, the ultrafilter $\alpha$ is uniquely determined by the sequence, and, having no minimal elements, constitutes its own almost-equality class in $\HH^\circ$. It follows that the class of $\alpha$ in $\Re\HH$ is maximal, and that, in particular, no class arising in the same manner may be smaller. We will presently show that the class $\alpha$ determines in $\Re\HH$ is, in fact, of codimension $1$, implying it is its own comparability component in $\rho(\bd X)$, and that its preimage under $\rho$ is a closed subset of $\bd X$. We shall call $\alpha$ a \emph{singular} point of $\HH^\circ$.
	\item[The sequence $(B_n)_{n=1}^\infty$ is eventually-constant.] In this case, denote the terminal value of the sequence by $B$, and recall every block $B$ is a closed convex subspace of $X$. By the definition of equivalence of chains, we may assume $B_n=B$ for all $n$. This allows us to consider the restriction $r^X_B:\HH\to\HH\big|_B$ together with the injective dual map of $\left(\HH\big|_B\right)^\circ$ into $\HH^\circ$. This map is a closed continuous embedding mapping almost-equality classes onto almost-equality classes, and we deduce that the almost-equality class of $\alpha$ cannot accumulate at a singular point of $\HH^\circ$. This  proves our previous assertion that the almost-equality class of a singular point of $\HH^\circ$ is of codimension $1$.

Now, since every element of $T(B)$ is transverse to every element of $P_c(B)$, we conclude that either $h_n\in T(B)$ for all $n$ or $h_n\in P_c(B)$ for all $n$. Since both $T(B)$ and $P_c(B)$ are the proper elements of poc-sets associated to trees, an ultrafilter on $\HH$ may not contain a pair of inequivalent chains of either type. Thus, the class corresponding to $\alpha$ in $\Re\HH$ is necessarily of codimension at most $2$, and we note that $\bd B$ is safely-connected with respect to $\HH\big|_B$. This implies that the safe path-components of any two block-boundaries are equal, showing that, except for the safe components corresponding to the ends of the Bass-Serre tree, there exists only one more safe path-component in $\bd X$. 
\end{description} 
To conclude this example, we see that the map $\rho$ does not distinguish among the boundaries of the spaces $X_\alpha$, though it retains essential information regarding the structure of all these spaces.

\subsection{Generalizing the examples} 
\subsubsection{Surjectivity of $\rho$.} In all the preceding examples we have witnessed situations where $\rho$ was essentially surjective: all non-principal classes were lying in the image of $\rho$. The following examples shows one cannot expect this to be true in general:
\begin{example}[$\rho$ is not necessarily surjective] Let $X$ be the subspace
\begin{equation}
	X=\bigcup_{n\in\NN}[n-1,n]\times[n^2,\infty),
\end{equation}
with the standard cubulation $\HH$ induced from $\EE^2$. There are three non-principal classes in $\Re\HH$. However, this space is not almost-extendible, and has only one boundary point corresponding to vertical rays -- both for the same reason: for any $C>0$, if $x$ is any point at a distance greater than $C$ to the vertical ray emanating from the point $x_0=(0,1)\in X$, then there is no geodesic ray $[x_0,\xi)$ passing through $B(x,C)$.
\end{example}
Note that by a theorem of Ontaneda \cite{[Ontaneda]}, a CAT(0) space admitting a geometric action by a group is almost extendible. This is why we are tempted to ask the following question:
\begin{question}\label{conj:rho is onto for cubings} Suppose $X$ is a proper cubing with standard halfspace system $\HH$. If $X$ is almost extendible, is it true that then every non-principal class in $\Re\HH$ lies in the image of the boundary representation map?
\end{question}
The results we have to report in this direction are much weaker, and serve, rather, as indications of situations in which one has an easy positive answer:
\begin{prop}\label{prop:codim 1 lies in the image of rho} Suppose $X$ is a proper cubing with standard halfspace system $\HH$, and let $\Pi\neq\Sigma\in\Re\HH$ be a class of codimension $1$. Then $\Sigma\in\rho(\bd X)$.
\end{prop}
\begin{thm}\label{thm:rho is onto for cubings} Suppose $X$ is a proper cubing with standard halfspace system $\HH$, and let $\Pi\neq\Sigma\in\Re\HH$. If, for some $\pi\in\Pi$ the canonical flow $\left(F^n(\pi)\right)_{n=1}^\infty$ fellow-travels a geodesic ray $[\pi,\xi)$ in $X$, then $\Sigma=\rho(\xi)$. 
\end{thm}
\proof{of prop \ref{prop:codim 1 lies in the image of rho}} Let $\Sigma\neq\Pi$ be an element of $\Re\HH$ and fix $\sigma_\infty\in\Sigma$. Let $(\sigma_n)_{n=1}^\infty$ be a sequence of ultrafilters in $\Pi$ converging on $\sigma_\infty$. Since $X$ is proper, considering the $\sigma_n$ as vertices of $X$ we may pass to a subsequence such that the $\sigma_n$ converge to a point $\xi\in\bd X$. As a result, $\sigma_\infty\in\HH^\xi$, and we have 
\begin{equation}\label{eqn:bad limits}
	\Pi<\rho(\xi)\leq\Sigma\,.
\end{equation}
In particular, $\rho(\xi)=\Sigma$ whenever $codim(\Sigma)=1$.\ep\\

\proof{of theorem \ref{thm:rho is onto for cubings}} We begin the proof in almost the same way as the previous one: let $\Sigma\neq\Pi$ be an element of $\Re\HH$, let $F=F_\Sigma$ denote the canonical flow of $\Sigma$ (defined on $\Pi$ -- see definition \ref{defn:canonical flow}), and let $\pi\in\Pi$ be an ultrafilter satisfying $F^n(\pi)\in N_R\left([\pi,\xi)\right)$ for all $n$ for some $\xi\in\bd X$. Let $\sigma_\infty=pr_\Sigma(\pi)=\lim_{n\to\infty}F^n(\pi)$, where the limit is taken in $\HH^\circ$ w.r.t. the Tychonoff topology.

As before, one must have $\rho(\xi)\leq\Sigma$. Suppose it were true that $\rho(\xi)<\Sigma$; then $P(\Sigma)\cap P(\xi)$ would have contained an infinite descending chain $(c_n)_{n=1}^\infty$. Without loss of generality we have $c_1^\ast\in\pi$, which, together with $c_1\in P(\xi)$ implies $[\pi,\xi)$ is (closure-)contained in $c_1^\ast$. However, for any $n>R$ we have $d(c_{n+1},c_1^\ast)>R$ (in $(X,d)$), which contradicts $F^{n+1}(\pi)$ lying in the $R-$neighbourhood of $[\pi,\xi)$.\ep
\begin{remark} Of course, the assumption of the last theorem need not regard approximating sequences (for $\Sigma$) generated by the corresponding canonical flow. However, the examples we have indicate that this is the correct formulation (consider the preceding example where $\rho$ is not surjective). The canonical flow seems to hold much information about the cubing $X$; for example, in any case when $\Delta(F_\Sigma(\pi),h)>\Delta(\pi,h)$ for all $h\in P(\Sigma)^\ast$ and $\pi\in\Pi$ (this is the case in $\EE^n$ with the standard cubulation), it is easy to show that the sequence $F^n(\pi)$ lies on a geodesic. In a way, this is our motivation for hoping that canonical flows do indeed fellow-travel geodesic rays -- at least for the case when $X$ admits a proper co-compact cellular action by a group. 
\end{remark}
\subsubsection{Giant components.}
The second feature common to the preceding examples was that the spaces and the cubings involved were all one-ended and the image of $\rho$ in all cases contained a unique `giant component' --
\begin{defn}[giant component] Suppose $\HH$ is a halfspace system on a proper CAT(0) space $X$. A comparability component $K$ of $\rho(\bd X)$ will be called a \emph{giant component}, if its preimage in $\bd X$ under $\rho$ is dense with respect to the cone topology.
\end{defn}
\begin{remark} As an example consider a situation when $\rho(\bd X)$ has a \emph{finite} giant component $K$ (like in the case of halfspace systems in $\EE^d$): in this situation, the closure formula implies $K=\rho(\bd X)$. The example of the standard presentation $2$-complex for $F_2\times\ZZ$ shows that $K=\rho(\bd X)$ is also possible for infinite giant components.
\end{remark}
\begin{thm} Suppose $\HH$ is a uniform halfspace system on a proper CAT(0) space $X$. If $\rho(\bd X)$ has a giant component, then $X$ is one-ended. 
\end{thm}
\begin{remark} The converse is not true: for example, we have seen that for $X=\FAT{H}^2$, any uniform halfspace system $\HH$ makes any two values of $\rho$ incomparable. Thus, every point in $\bd X$ corresponds to its own comparability component in spite of $X$ being one-ended.
\end{remark}
\proof{} Let $K$ be a comparability component of $\rho(\bd X)$, and let $Y=\cl{\rho^{-1}(K)}$, the closure taken with respect to the cone topology. We show that any two points of $\xi,\eta\in Y$ belong to the same end of $X$: taking $U$ and $V$ to be cone neighbourhoods of $\xi$ and $\eta$ respectively, both not intersecting a ball $B(x_0,r)\subset X$ of a given radius $r>0$, we find points $\xi'\in\rho^{-1}(K)\cap U$ and $\eta'\in\rho^{-1}(K)\cap V$. For some $R>r$ it is then possible to connect $\xi'(R)$ with $\eta'(R)$ with a piecewise-circular path. This path can be then extended on both ends by rays asymptotic to $[x_0,\xi)$ and $[x_0,\eta)$, which, by properties of rays in CAT(0) spaces, will necessarily be contained in $U$ and $V$, respectively.

In particular, if $K$ is a giant component, $X$ is one-ended.\ep
\begin{thm}[uniqueness of giant components]\label{thm:uniqueness of giant components} Suppose $\HH$ is a uniform halfspace system in a proper CAT(0) space $X$, and suppose $F$ is a closed convex subspace of $X$ such that
\begin{enumerate}
	\item $F$ coarsely separates $X$, and
	\item $\rho_F\left(\bd F\right)$ has a unique comparability component.
\end{enumerate}
Then $\rho(\bd X)$ has at most one giant component.
\end{thm}
\proof{} The first assumption on $F$ means that $\bd F$ separates $\bd X$. It will be enough to show that, if $K$ is a giant component of $\rho(\bd X)$ then $\rho(\bd F)$ is contained in $K$. Since different comparability components of $\rho(\bd X)$ do not intersect, it will follow that $K$ is the only giant component in stock.

Now, the second assumption on $\rho(\bd F)$, together with proposition \ref{prop:the restriction equality}, imply it is enough to show that $\rho(\bd F)$ intersects $K$.

Now, since $F$ coarsely separates $X$, $\bd F$ coarsely separates $\bd X$. Since $\bd F$ is closed in $\bd X$, we may write $\bd X\minus\bd F=U\cup V$, where $U$ and $V$ are disjoint open subsets of $\bd X$. Since $K$ is a giant component, there are points 
\begin{equation}
	\alpha\in \rho^{-1}(K)\cap U\,,\quad\beta\in\rho^{-1}(K)\cap V\,
\end{equation}
and $\bd F$ intersects any (safe) path from $\alpha$ to $\beta$ in $\bd X$. In particular, $\rho(\bd F)$ intersects $K$, and we are done.\ep
\begin{cor} Suppose a group $G$ acts geometrically on a CAT(0) space, and suppose $X$ admits a uniform $G$-invariant halfspace system $\HH$. If $G$ has a codimension-one free-abelian subgroup $A$ of rank at least $2$, then $\rho(\bd X)$ has at most one giant component, which, when it exists, is characterized as the comparability component containing the $\rho$-images of the boundaries of all coarsely-separating flats in $X$.
\end{cor}
\proof{} By the flat torus theorem, there exists a flat $F$ of dimension $r\geq 2$ in $X$, on which $A$ acts co-compactly by translations; since $A$ is a codimension-one subgroup, $F$ coarsely separates $X$. Proposition \ref{prop:Euclidean boundaries are safe} then implies that $\bd F$ is safely connected, allowing to apply the last theorem and its proof.\ep

\section{Co-compact actions.}
Suppose, once again, that $G$ is a group acting geometrically on the CAT(0) space $X$, and $X$ has a $G$-invariant uniform halfspace system $\HH$. In this section our aim will be to use the decomposition mapping $\rho:\bd X\to\Re\HH$ for establishing when $G$ acts co-compactly on the cubing dual to $\HH$.
\begin{thm}\label{thm:co-compactness criterion} Suppose $G$ is a group acting geometrically on a CAT(0) space $X$, and suppose $\HH$ is a $G$-invariant uniform halfspace system. If every non-principal class of $\Re\HH$ lies in the image of the boundary decomposition map $\rho:\bd X\to\Re\HH$, then the action of $G$ on $C(\HH)$ is co-compact.
\end{thm}
Note that the converse to this result is, in fact, a harder version of the problem discussed in the preceding section (question \ref{conj:rho is onto for cubings}, theorem \ref{thm:rho is onto for cubings}).\\

The first step in our discussion of the problem will be to reduce it to a combinatorial problem, which is done in a way analogous to that used by Williams in his thesis \cite{[Wil]}. Most of the technical details of this reduction, as well as some of the required results, were dealt with by the author in \cite{[Gu-loc-fin]}. The most important of these results states that $C(\HH)$ is locally finite.\\

In the following paragraph we provide an overview of the results and notions we shall use from that source.

\subsection{Consistent ultrafilters.}
For any point $x\in X$ let us consider the set $B_x$ of all $h\in\HH$ satisfying $x\in h$. $B_x$ is a filter base. We define:
\begin{defn}\label{defn:consistent set} A point $x\in X$ is said to \emph{support} a subset $A\subseteq\HH$, if $x\in\cl{h}$ for all $h\in A$. A subset $A\subset\HH$ is said to be \emph{consistent}, if it has a supporting point. We shall adopt the convention that the empty set is inconsistent.
\end{defn}
For ultrafilters, consistency is an easy matter:
\begin{lemma}[consistent ultrafilters] Suppose $\HH$ is a halfspace system in a CAT(0) space $X$. Then:
\begin{enumerate}
    \item any consistent ultrafilter is principal;
    \item any point $x\in X$ supports an ultrafilter.
\end{enumerate}
The set of all consistent ultrafilters will be henceforth denoted by $\Pi_0$.\ep
\end{lemma}
We see that the consistent ultrafilters -- or, at least, those which are of the form $B_x$ for some $x\in X$ -- correspond to chambers. This, in effect, is what provides us with the combinatorial reduction of the co-compactness problem.
\begin{lemma}\label{lemma:consistent principal UF's are G-finite}
The set $\Pi_0$ is $G$-finite. In particular, there exists a constant $D>0$, such that if $\pi,\pi'\in\Pi_0$ are supported on the same point, then $|\pi\vartriangle\pi'|\leq D$.\ep
\end{lemma}
Recall the metric $\Delta$ on $\Pi$. It is $G$-invariant. Now, since $G$ acts on $X$ stabilizing $\Pi_0$, it also stabilizes the level sets $\Pi_n$ ($n\in\NN\cup\{0\}$) of the function $\Delta(-,\Pi_0):\sigma\mapsto \Delta(\sigma,\Pi_0)$. Since the action of $G$ on $\Pi_0$ is co-finite, the action of $G$ on $\Pi$ will be co-bounded if and only if the function $\Delta(-,\Pi_0)$ is bounded. Since $C(\HH)$ is locally finite, its quotient by $G$ will be compact if and only if it is bounded. Therefore, studying the growth of $\Delta(-,\Pi_0)$ will be our tool for studying the co-compactness problem.\\

At the base of the technique lies a geometric interpretation of $\Delta(\pi,\Pi_0)$:
\begin{lemma}\label{lemma:characterizing level sets} Suppose $\pi\in\Pi$. Then $\Delta(\pi,\Pi_0)\leq n$ if and only if there exists a subset $A$ of $\pi$ of size $n$ such that $\pi\minus A$ is consistent.\ep
\end{lemma}
Then, one needs a tool for understanding how $\Delta(\pi,\Pi_0)$ changes as one performs walks on the $1$-skeleton of $\Pi$. For this we define --
\begin{defn} Suppose $\pi\in\Pi$ and $a\in\min(\pi)$. Denote
\begin{eqnarray*}
	a\in\min(\pi)_+ &\IFF& \Delta([\pi]_a,\Pi_0)>\Delta(\pi,\Pi_0)\,,\\
	a\in\min(\pi)_- &\IFF& \Delta([\pi]_a,\Pi_0)<\Delta(\pi,\Pi_0)\,.
\end{eqnarray*}
Note that $\pi\in\Pi_0$ iff $\min(\pi)_-$ is empty.\ep
\end{defn}
The basic observations regarding the signed minimal sets of $\pi\in\Pi$ are:
\begin{lemma}\label{lemma:min pi minus is inconsistent} For all $\pi\in\Pi$, the set $\min(\pi)_-$ is inconsistent.\ep
\end{lemma}
\begin{cor}[three ways to go down]\label{cor:three ways to go down} If $\pi\notin\Pi_0$, then $\min(\pi)_-$ contains at least three distinct elements.\ep
\end{cor}
An important means of assessing distances to $\Pi_0$ is the growth of the following objects (as $\pi$ recedes from $\Pi_0$):
\begin{defn}[Shadows]\label{defn:shadow} For all $\pi\in\Pi$, let the \emph{shadow of $\pi$} be defined as the set
\begin{equation}
	\sh{\pi}=\left\{\sigma\in\Pi_0\,\big|\,\delta(\sigma,\pi)=\Delta(\sigma,\Pi_0)\right\},
\end{equation} 
and let the \emph{dual shadow} $\shcirc{\pi}$ of $\pi$ be defined to be
\begin{equation}
	\shcirc{\pi}=\left\{h\in\HH\,\big|\,\sh{\pi}\subseteq S_h\right\}.
\end{equation}
\end{defn}
Observe that $\shcirc{\pi}$ is a filter-base, and it is natural to expect that $\shcirc{\pi}$ be contained in $\pi$. If that is the case, it would mean that, as the distance of $\pi$ from $\Pi_0$ increases, $\shcirc{\pi}$ diminishes accordingly, testifying to the growth of $\sh{\pi}$. This is exactly what we will require for the study of the co-compactness problem.\\

In order to gain a feeling of why shadows should be related to the co-compactness problem, let us consider once again the example of a Coxeter group $(W,R)$ acting on its Davis-Moussong complex (and the associated halfspace system -- call it $\HH$ -- which is, in fact, isomorphic to the corresponding system, ordered by the domination order introduced in \cite{[BriHow]} for positive roots, and extended to all roots by \cite{[NibRee1]}). It is now known, through the work of \cite{[Wil]} and \cite{[Caprace]}, that the action of $W$ on $C(\HH)$ is co-compact if and only if $W$ does not contain a Euclidean triangle subgroup. Let us consider such a subgroup $W'$, and its associated Davis-Moussong (sub)complex, which is isometric to a Euclidean plane $P$ tessellated by hexagons (corresponding to the finite dihedral parabolic subgroups): one notices that an inconsistent ultrafilter $\pi$ will necessarily have $\min(\pi)_-$ define (through intersection) a triangle in $P$, and that these triangles can grow arbitrarily large; the fact that these triangles come in infinitely many different sizes shows that the number of conjugacy classes (in $W'$, and hence in $W$) of subgroups $W''<W$ isomorphic to $W'$ must be infinite (as $W$ acts on $X$ by isometries) -- which hints at Wilson's co-compactness criterion ($W$ acts co-compactly on $C(\HH)$ iff $W$ contains only finitely many conjugacy classes of triangular subgroup of each admissible type). Wilson and Caprace went on to study the structure of relations among reflections creating this abundance of `similar shapes' (The hyperbolic triangle groups, lack this abundance of similar shapes, for the obvious geometric reasons; in fact, Caprace's work exhausts the possibilities for such `shapes' in the general case, showing that the only reason for non-co-compactness is the presence of infinite similarity classes of Euclidean `shapes'). In our, more general case, however, we shall have to study the `shapes' themselves -- these are modelled by shadows.

\begin{figure}[htb]
	\includegraphics[width=\textwidth]{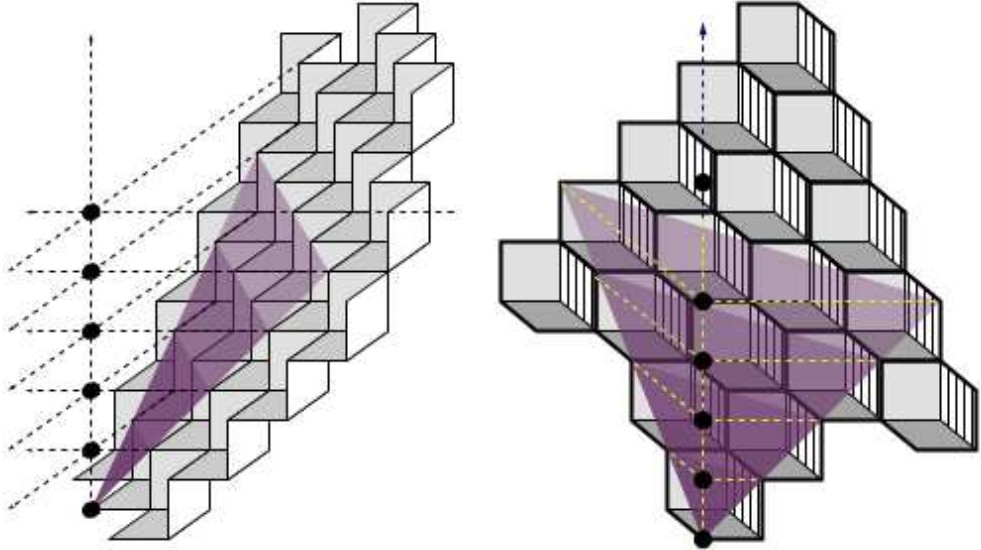}
	\caption{\protect\scriptsize The hexagonal packing in $\EE^2$ is the Davis-Moussong complex of a Coxeter triangle group for which $C(\HH)$ is isomorphic to the standard cubulation of $\EE^3$. Each vertex $v$ of the packing corresponds to the ultrafilter of halfspaces containing $v$; so that the Cayley graph of $W$ (which is the $1$-skeleton of the packing) is embedded in $C(\HH)\cong\EE^3$ as shown above (thick lines, right view), while inconsistent ultrafilters give rise to shadows.\protect\normalsize}
	\label{figure:shadows}
\end{figure}

The technicalities we require are summarized in a result we had proved in \cite{[Gu-loc-fin]}:
\begin{lemma}\label{lemma:bounds on shadows} For every $\pi\in\Pi$ one has
$\min(\pi)_+\subseteq\shcirc{\pi}\subseteq\pi$. Furthermore, for any $a\in\min(\pi)$ one has:
\begin{enumerate}
	\item if $a\in\min(\pi)_+$ then $\sh{\pi}\subsetneq\sh{[\pi]_a}$ and $\shcirc{[\pi]_a}\subsetneq\shcirc{\pi}$;
	\item if $a\in\min(\pi)_-$ then $\shcirc{\pi}\cap\min(\pi)_-=\varnothing$.
\end{enumerate}
\end{lemma}

\subsection{Proof of theorem \ref{thm:co-compactness criterion}}
\subsubsection{Limits of geodesic rays.}
We need a geometric lemma involving the structure of $T(\xi)$ in a uniform halfspace system.
\begin{lemma} Suppose $A\subseteq T(\xi)$ is a non-empty set containing elements $a_1,\ldots,a_m$ such that $A\minus\{a_1,\ldots,a_m\}$ is consistent. Then $A$ is a consistent set.
\end{lemma}
\proof{} It is enough to prove the result for the case $m=1$. Suppose $A$ were inconsistent. This would mean that 
\begin{equation}
	Y:=\bigcap_{a\in A\minus\{a_1\}}\cl{a}\subset a_1^\ast.
\end{equation}
Now, on one hand, $Y$ is consistent, closed and convex, and so it contains a geodesic ray $[y,\xi)$, which, therefore lies entirely in $a_1^\ast$. On the other hand we see that $\cl{a_1}$ contains a ray $[x,\xi)$, implying $a_1\in P(\xi)$ -- contradiction.\ep\\

\begin{figure}[htb]
	\includegraphics[width=\textwidth]{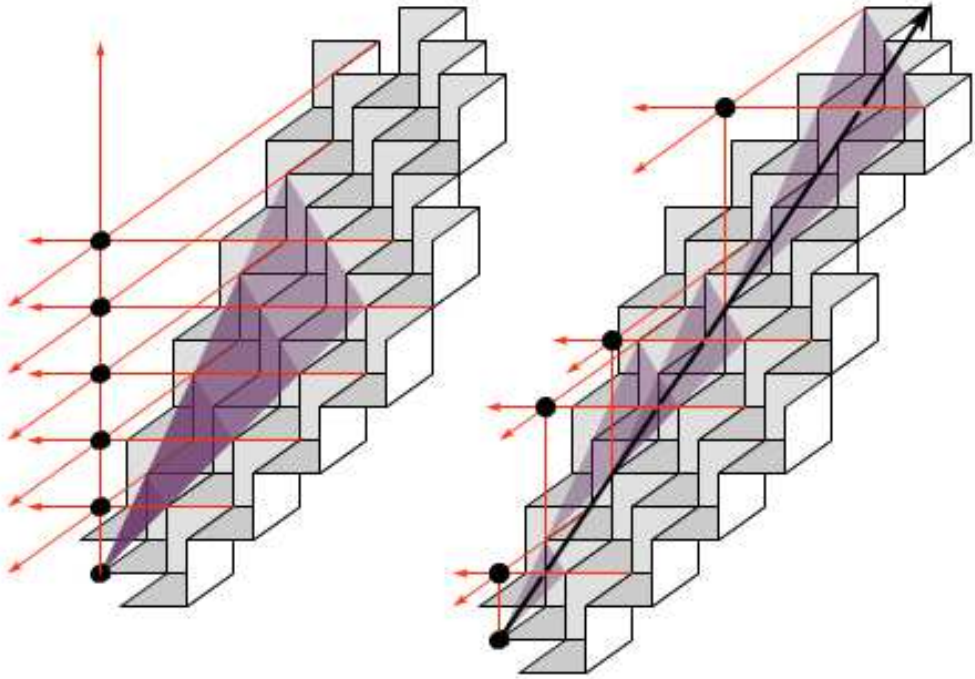}
	\caption{\protect\scriptsize The hexagonal packing again. Shadows of ultrafilters receding to an unbounded distance from $\Pi_0$ (left) will increase without moving away from some common point. Shadows of ultrafilters converging (tamely) to $\rho(\xi)$, $\xi\in\bd X$ (right) will `move faster than they grow'.\protect\normalsize}
	\label{figure:escaping shadows}
\end{figure}

The following proposition is the crucial ingredient in the proof of theorem \ref{thm:co-compactness criterion}. In a way, it states which geodesic vertex paths in the one-skeleton of the cubing dual to $\HH$ produce limits representing elements of $\Re\HH$ which do not lie in the image of $\rho$. The motivation for this proposition was given by the example of the hexagonal packing, as shown in figure \ref{figure:escaping shadows}.
\begin{prop}[escaping vertex paths]\label{prop:escaping geodesics}
Suppose $\Pi$ contains a geodesic sequence $(\pi_n)_{n=0}^\infty$ satisfying
\begin{displaymath}
	\Delta(\pi_{n+1},\Pi_0)>\Delta(\pi_n,\Pi_0)
\end{displaymath}
for all $n$, and let $\Sigma$ denote the almost-equality class of the limit $\pi_\infty$ of this ray. Then $\Sigma\notin\rho(\bd X)$.
\end{prop}
\proof{} Suppose, on the contrary, that $\Sigma=\rho(\xi)$ for some $\xi\in\bd X$. 

By lemma \ref{lemma:geodesic sequences converge} and corollary \ref{cor:geodesic sequences converge tamely}, truncating an initial segment of the sequence $(\pi_n)_{n=0}^\infty$ we may assume that $\pi_n\cap P(\xi)$ is constant (and then, for all $n$, we would have $\pi_n\cap P(\xi)=\pi_\infty\cap P(\xi)$, which is a principal ultrafilter on $P(\xi)$), implying that all the elements of $\HH$ whose orientations are being reversed along this ray of $\Pi$ are halfspaces belonging to $T(\xi)^\ast$. Let us denote these halfspaces by $a_n$, so that we have $\pi_{n+1}=\left[\pi_n\right]_{a_n}$ and $a_n^\ast\in\min(\pi_{n+1})_-$ for all $n$.

Now observe that, for any $n$ and $h\in\pi_n$, if $h\notin\shcirc{\pi_n}$, then $h\in\pi_{n+1}$ (otherwise, this would imply that $h\in\min(\pi_n)_+$, producing $h\in\shcirc{\pi_n}$). However, by lemma \ref{lemma:bounds on shadows}, the set $\shcirc{\pi_{n+1}}$ is contained in $\shcirc{\pi_n}$, implying $h\notin\shcirc{\pi_{n+1}}$. Applying induction we see that $h$ then belongs to $\pi_{n+k}$ for all $k$.

As a result, for all $n$ we must have $\pi_n\minus\shcirc{\pi_n}\subseteq T(\xi)$. In particular, $\min(\pi_n)_-$ is contained in $T(\xi)$ for all $n$, as a result of applying corollary \ref{lemma:bounds on shadows}.

Let $\pi=\pi_1$, and recall that $\min(\pi)_-$ is inconsistent (lemma \ref{lemma:min pi minus is inconsistent}). This contradicts the preceding lemma when applied with $A=\min(\pi)_-$.\ep

\subsubsection{Finishing the proof.}
\proof{of theorem \ref{thm:co-compactness criterion}}
Suppose $\Delta(-,\Pi_0)$ is unbounded, and let us construct a geodesic ray $(\pi_n)_{n=1}^\infty$ such that $\Delta(\pi_{n+1},\Pi_0)>\Delta(\pi_n,\Pi_0)$ for all $n$. Applying proposition \ref{prop:escaping geodesics} will produce a point of $\Re\HH$ outside the image of $\rho$. Our construction of the desired ray in $\Pi$ is a pigeonhole argument slightly generalizing an argument Williams had used for the special case of Coxeter groups \cite{[Wil]}.

Since $\Delta(-,\Pi_0)$ is unbounded, for every $n\in\NN$ we select an ultrafilter $\pi_n$ at a distance $n$ from $\Pi_0$. For each $\pi_n$, let $\pi_{0,n}$ be an element of $\Pi_0$ at a distance $n$ from $\pi_n$. Since $G$ acts co-finitely on $\Pi_0$ and preserves $\Delta(-,\Pi_0)$, we may select a subsequence $(\pi_{n_k})_{k=1}^\infty$ such that $\pi_{0,n_k}$ is the same for all $k$. Denote this ultrafilter by $\pi_0$. 

Since $\Delta(\pi_{n_k},\pi_0)=\Delta(\pi_{n_k},\Pi_0)=n_k$, for each $k$ there is a vertex-path $p_k$ in $\Pi$ of the form 
\begin{equation}
    p_k=(\pi_0,\pi_{1,n_k},\ldots,\pi_{n_k-1,n_k},\pi_{n_k})
\end{equation}
satisfying $\Delta(\pi_{i,n_k},\Pi_0)=i$ for all relevant $i$.

By the local finiteness of $C(\HH)$, one may pass to a subsequence yet again, so that every path $p_{k_l}$ is an initial subpath of the path $p_{k_{l+1}}$. The union of all the $p_{k_l}$ is then an escaping geodesic ray in $\Pi$.\ep

\bibliographystyle{alpha}
\bibliography{prop}

\end{document}